\documentclass[11pt]{article}
\usepackage{euscript}
\usepackage{amsfonts}
\usepackage{sidecap}
\usepackage{hyperref}
\usepackage{url}
\usepackage{mathptmx}
\usepackage{times}
\usepackage{amsmath,amssymb,amsthm,euscript,latexsym}
\usepackage{mathrsfs}
\usepackage{bm}
\usepackage[dvips]{graphicx}

\usepackage{exscale}
\usepackage{color}
\usepackage{tikz}
\usetikzlibrary{calc,patterns,decorations.pathmorphing,decorations.markings}
\usetikzlibrary{snakes,arrows,shapes}
\usepackage{pgf}
\graphicspath{{images/}}

\newcommand{\bmat}{\left[ \begin{matrix}}
\newcommand{\emat}{\end{matrix} \right]}

\newcommand{\beq}{\begin{equation}}
\newcommand{\eeq}{\end{equation}}

\newcommand{\nd}{\noindent}
\newcommand{\q}{\quad}
\newcommand{\qq}{\qquad}

\newcommand{\Frac}[2]{{\displaystyle\frac{#1}{#2}}}
\newcommand{\met}{\mathop{\frac{1}{2}}\,}
\newcommand{\Span} {\mbox{\rm span}\,}

\newcommand{\E} {\mathbb{E}\,}
\newcommand{\Cbb} {\mathbb{C}}
\newcommand{\Rbb} {\mathbb{R}}

\newcommand{\Zbb} {\mathbb{Z}}

\def\v{{\rm v}}
\def\ib{{\mathbf i}}

\def\bb{{\mathbf b}}
\def\ab{{\mathbf a}}

\def\cb{{\mathbf c}}

 \def\qb{{\mathbf q}}
\def\pb{{\mathbf p}}

\def\xb{{\mathbf x}}
\def\vb{{\mathbf v}}
\def\yb{{\mathbf y}}

\def\Eb{{\mathbf E}}
\def\Hb{{\mathbf H}}

\def\Xb{{\mathbf X}}

\def\Sb{{\mathbf S}}

\def\xib{\boldsymbol{\xi}}

\def\taub{\boldsymbol{\tau}}

\def\varphib{\boldsymbol{\varphi}}

\newtheorem{thm}{\bf Theorem}

\newtheorem{defn}{\bf Definition}
\newtheorem{lem}{\bf Lemma}
\newtheorem{cor}{\bf Corollary}

\newtheorem{rem}{\bf Remark}

\begin{document}

\thispagestyle{empty}

\begin{center} 
{\textcolor{blue}{\Huge {  On Irreversibility\\[4. mm] and Stochastic Systems\\ [4. mm] Part One}}} \\[1.cm]

 in memoriam of  Jan. C. Willems\\[1.cm]

{\Large   GIORGIO PICCI}\\[.6cm]
Department of Information Engineering\\
  University of Padova, Italy\\
  \tt{picci@dei.unipd.it}

\vspace{2. cm}
\end{center}

\nd{\bf Abstract:} We attempt to characterize irreversibility of a dynamical system from the existence of different forward and backward mathematical representations depending on the direction of the time arrow. Such different representations have been studied intensively and are shown to exist  for stochastic diffusion models. In this setting one has however to face the preliminary justification of stochastic description for physical systems which are described by classical mechanics as inherently deterministic and conservative.\\
 In part one of this paper we first address this modeling problem  for linear systems in a deterministic context. We show that forward-backward representations can also describe conservative finite dimensional deterministic systems when they are  coupled to an infinite-dimensional conservative heat bath. A novel key observation is  that  the heat bath acts on the finite-dimensional conservative system by {\em state-feedback}  and can shift its eigenvalues to make the system dissipative but  may also generate another  totally unstable model which naturally evolves backward in time. 
 
 In the second part, we  address the stochastic description of these two representations. Under a natural family of invariant measures the heat bath can be shown to  induce  a white noise input acting on the system making it look like a true dissipative diffusion.

\section{Introduction}

Irreversibility is a physical phenomenon which so far  has escaped a precise definition \cite{Schrodinger-51}, \cite{Lebowitz-93-99}, \cite{Lebowitz-08}. It has to do with the impossibility of describing by a unified dynamical law the motion of a system both in the forward and in the {\em backward} direction of time, say of a particle of smoke exiting from a cigarette,  describing how that particle should travel back into the cigarette tip when the direction of time is reversed.  One could actually venture to say that {\em most} physical systems   behave very differently when the time arrow is reversed. In this sense, we may agree to say that they are not reversible.  \\
In this essay, inspired by Galileo's tenet that nature can only be described faithfully by mathematics, we shall only attempt to    propose a  {\em mathematical} definition of irreversibility. This definition  will be stated in the language of  probability and on the theory of probabilistic (also called stochastic)   modeling of physical systems as it  has been developing in theoretical Systems Engineering in recent years, as summarized  e.g. in the book \cite{LPBook}. A  detailed physical analysis or justification of the link of the ideas proposed in this report with observed irreversible phenomena will be left a bit in the vague and delegated to the expertise and intuition of the reader.\\
   For reasons of simplicity and clarity we shall mostly  restrict to stationary  linear phenomena although a generalization to transient and nonlinear systems is likely to be possible, based on some general first  principles that we shall   try to unveil by discussing some simple examples. These principles are stated in the setting and language of   Systems Theory  which may sound extraneous to some physicist although in my opinion this setting  turns out to be very natural in helping to streamline the basic conceptual background of the phenomenon. The basic  underlying issue turns actually out to be that, in order to explain   irreversibility  one needs to face the preliminary question of how from an intrinsically deterministic conservative system like those of classical mechanics  (conservative and hence {\em reversible}) one can  extract an, at least  local, dissipative (and hence {\em irreversible}) finite-dimensional {\em stochastic} dynamical model.

It is by now well-known that  stochastic diffusion models  entail  a basic duality, in fact a fundamental structural difference,  between   two possible dynamic representations of the {\em same  stochastic process} when the time arrow is reversed \cite{Nelson-58,Lindquist-P-79,Lindquist-P-85b,Protter-87}. Although they  describe the {\em same} stochastic system, there is a fundamental diversity of  the  two possible dynamic forward/backward  representations which  arise when the time arrow is reversed.  The phenomenon   generalizes to time-varying and nonlinear models and we wan to suggest that it should  be considered as a general  {\em mathematical} elucidation of the phenomenon of irreversibility. Perhaps this may seem a somewhat too  drastic  conceptual simplification of a rather mysterious and intricate physical fact. To convince the reader one should ideally arrive to relate logically the empirical idea of irreversibility in classical thermodynamics to this rather sophisticated mathematical phenomenon through a consistent sequence of arguments. We shall act with this goal in mind. During the course of this program we shall however have to discuss  also other important intermediate conceptual steps.

A warning and an apology are in order regarding the somewhat simplistic and brute-force derivation of some key facts regarding statistical mechanics  and  the role of Entropy in the historical introduction section. We hope nevertheless that this naive brute force approach will make the motivations  of this work to emerge with  clarity.

\subsection{ \bf On the Role of  Probability and Statistics}

\d{\em Background:} Probability theory is based on a mathematical conceptual model of an experiment, a {\em Probability Space}: $\{\Omega,{\mathcal A}, P\}$ where $\Omega$ is an ideal ``urn" of states $\omega\in\Omega$ of the ``environment"  chosen by nature in each experiment.

\nd Jan C. Willems, a very respected thinker,  in {\em Reflections on Fourteen Cryptic Issues Concerning the Nature of Statistical Inference} \cite{Willems-03}  writes:

{\em In industrial and econometric systems there is no urn; even an hypothetical repetition (of the experiment) is  impossible...Thus 
every physical system is {\bf ``deterministic !!"}; probability can just be used as an approximation device.}\footnote{May here quote Einstein's famous statement about God not playing dice where however,  he refers to the axiomatic structure  of quantum theory which is extraneous to and far beyond  the classical thermodynamical setting we are interested in here. }

\medskip

{\bf However:} there is huge evidence that every physical system shows an {\em irreversible behaviour}; examples are everywhere; see e.g. \cite{Lebowitz-08}. To extract deterministic mathematical laws from  experiments you need to artificially clean out  the data pretending there is no drag, no disturbances etc..  This was first done by Galileo for determining the square law of a falling body. Because of drag and disturbances {\em no physical system can behave the same in reverse time}. Instead, the standard deterministic models like ordinary differential equations are {\em reversible}, for them time can run both ways, but {\em real physical systems are not}.

Any attempt to reach some understanding of irreversibility starts from the notion of {\bf Entropy}. We shall retrace the basic idea in a very streamlined way in the next chapter.

\section{\bf Historical Introduction:  Entropy}
The  concept of entropy,   introduced by  {\bf Rudolf Clausius} in 1865 \cite{Clausius-865} has been object of a very long debate and equivalent reformulations see e.g. \cite{Bernstein-011}. Here we take a somewhat simplified point of view.

Let $ \Frac{q}{T}$ be the relative heat flow  supplied to a system at temperature $T$ in a transformation from a thermodynamic ({\em macroscopic}) initial  state $x_0$ (say $(p_0,V_0)$ for an ideal gas) to reach a state $x$. 
\begin{defn}\label{DefEntropy}
If the transformation  $x_0\to x$ is {\bf reversible}\footnote{Reversibility  in Physics is defined rather loosely. We shall mean that the state trajectory is generated by a deterministic smooth lossless dynamical system.}, the line integral in the phase space, along any path connecting $x_0$ to $x$ is independent of the path  and then defines a   function $S(x)$ of the state $x$,  the {\bf Entropy} of the system, by the relation:
\beq \label{EntrDef}
S(x)-S(x_0)= \int_{x_0}^{x}\Frac{q}{T}\,\mu(dx)\,.
\eeq 
\end{defn}
Hence by identifying the heat flow $q$ with the input $u$ which causes the temperature $T$ of the underlying (thermodynamical) deterministic dynamical system to change, entropy is seen as   a {\em storage function} in the sense of  J C Willems \cite{Willems-72}, a mathematically well-defined  object.

\begin{rem}Note that the function $S(x)$ is defined for {\em all reachable states of the system} which could be reached by a  ``reversible" transformation  or not;  it is unique up to an additive constant and is independent of the path followed  to reach the state $x$ .
\end{rem}

N.B. The condition for existence of the function $S(x)$  in Definition \ref{DefEntropy} in the language of System Theory  it is a condition  equivalent to {\em Lossless cyclo-dissipativity} i.e.
 $$
 \oint  \Frac{q}{T}\,\mu(dx) =0
$$
for every closed path in the phase space \cite{Willems-72}.
 The  classical (macroscopic) proof is based on Carnot Cycle.

Note  the key property that Entropy is  a {\bf function of the state} say $x=(p,V)$ for an ideal gas and hence makes perfect sense even if you can reach that state $x$ in an arbitrary way.

However in real physical systems all transformations entail a dissipative use of energy, resulting in inherent loss of usable heat when work is done, since  there is always  heat produced by internal  friction.
 
 \nd{\bf The second law of Thermodynamics} states that   for all physical  transformations 
 $$
 S(x)-S(x_0)\geq\,\int_{x_0}^{x}\Frac{q}{T}\,\mu(dx)
 $$
 that is, {\bf  Entropy is always a non-decreasing function of the state},  in particular  if the system is isolated ($q=0$) then $S(x)\geq\,S(x_0)$. Only if the transformation is reversible   the entropy stays constant, i.e. $S(x)=\,S(x_0)$. However for {\em real}   thermodynamic transformations  there is always loss of internal energy (heat) due to friction and we {\em always have strict inequality} i.e.  $S(x)>\,S(x_0)$. This is often taken as a {\em mathematical characterization  of irreversibility}. If (and only if) the transformation is {\bf irreversible} then Entropy is strictly increasing.

The {\em Mechanical Explanation of Irreversible Processes} is a century long unsolved problem which has generated a harsh debate among physicists  in the late 1800's. See the reference \cite{Steckline-83}.

General fact: friction or electrical resistance (i.e. physical irreversibility) is always described as  a {\em macroscopic} phenomenon. Can one explain irreversibility from its origin  as a {\bf microscopic phenomenon}?

Basic Difficulty:  Ideal gases of particles or electrons moving in an empty space cannot exercise friction. They must obey Newton's laws {\bf exactly} and their motion is by their very  nature {\em reversible}.

\bigskip

\subsection{Boltzmann's Program:}
 {\bf Probability can explain irreversibility!}

Starting point: The famous {\em  Maxwell–Boltzmann distribution}    a cornerstone of the kinetic theory of gases,  provides an   explanation of many fundamental macroscopic gaseous properties, including pressure and diffusion. The distribution was first derived by Maxwell in 1860 on heuristic grounds. Ludwig  Boltzmann later, in the 1870s, carried out precise investigations into the physical origins of this distribution. See \cite{Boltzmann-895}.\\
The Maxwell–Boltzmann distribution applies to rarified ideal gas particles; it describes the  {\em magnitude $\v$ of the velocity}  of the particles. Collisions are extremely rare and can be neglected.\\
 The basic assumption is that each component of the velocity vector    of any particle selected at random in   the ensemble  can be described   as   random variables having a {\bf Gaussian distribution}. 

Precisely: the Maxwell-Bolzmann distribution  is based on    the assumption that the three Cartesian components $ \v_1, \v_2, \v_3$ of the (vector) velocity of a particle are independent random variables each having the  same  Gaussian distribution
$$
p( \v_i)= (2\pi\, kT/m)^{-1/2}\exp{-\Frac{m  \v_i^2}{2kT}} , \qquad i=1,2,3
$$
where $\Frac {kT}{m}$ has the meaning of a   variance $(\equiv \sigma^2)$. The introduction of the term $kT$ in the variance  is due to Boltzmann.\\
Hence $\v_i/\sigma \sim {\mathcal N}(0,1)$   and the normalized square speed has a chi-squared istribution
$$
\frac{\v^2}{\sigma^2}:=  \sum_{i=1}^3 (\v_i/\sigma)^2 \sim \chi^2(3)\,.
$$
 By the change of variable $ x=\frac{\v^2}{\sigma^2}$ in the expression of $\chi^2(3)$ we obtain
$$
p(\v^2)= 4\pi(m/2\pi kT)^{3/2}v^2\exp{-\Frac{m v^2}{2kT}} 
$$
which is exactly the {\em Maxwell-Boltzmann distribution}.

  Note that the kynetic energy of a particle is 
 $$ 
 \Eb:=\met m\v^2=\met kT\frac{\v^2}{\sigma^2}\q \sim \q \met kT \,\chi^2(3)
 $$  
which agrees with Physics. In fact, from the mean of the $\chi^2(3)$ distribution (= the number of degrees of freedom) we immediately get the  relation
$$
\E\met m\v^2 = \frac{3}{2}kT.
$$
which is a well-known experimental relation of the average kinetic energy per molecule of (monoatomic) gases.

\subsection{\bf Equilibrium and  The velocity process}
 
Each particle of an isolated ideal gas system of $N$ particles must obey the laws of classical mechanics. The dynamics of such  system is generated by a quadratic Hamiltonian
$$
H(q,p)= \met \sum_{k=1}^N\|p_k\|^2 + \met \sum_{k,h=1}^Nq_k^{\top} V^2_{k,h}q_h
$$
where the potential matrix $V^2:= \left [V^2_{k,h}  \in \Rbb^{3\times 3} \right]_{k,h = 1,\ldots,N }$ is assumed symmetric positive definite.  The canonical equations are:
$$
\begin{cases}
\dot q_k(t) & = p_k(t)\\
\dot p_k(t) &= -  \sum_{h} V^2_{k,h}\,q_{h}(t)\,,\quad k\in \Zbb 
\end{cases}
\quad \Leftrightarrow \begin{cases}
\bmat \dot q(t)\\ \dot p(t)\emat= \bmat 0  & I\\ - V^2 & 0\emat \,\bmat q(t) \\p(t) \emat
\end{cases}
$$ 
where the components  $q_k, p_k$ of the canonical variables are the 3-dimensional space coordinate vectors and the relative momenta.
Each trajectory of the system must obey these equations and is uniquely defined given some initial condition say at time $t=0$.

The solution, called the {\em Hamiltonian flow}, $\Phi(t); t\in \Rbb$ is a group of linear operators on   the phase space which preserves the total energy $H(q,p)$. Any $\Phi(t)$-invariant probability measure on the phase space defines a situation of equilibrium.
 \begin{thm}
Any absolutely continuous $\Phi (t)$-invariant probability measure is a member of the  one-parameter family of densities
$$
\rho_{\beta} (q,p) =C\,exp[-\frac{1}{\beta}H(q,p)] ,\qq \beta>0 
$$
 where $C$ is a normalization constant and the parameter $\beta$ has been linked by Boltzmann to  the absolute temperature by the relation $ \beta:= kT $.
\end{thm}
\begin{proof}
We report a  proof  based  on Liouville theorem, see e.g. \cite[p. 428]{Goldstein-80}.\\
Since $\rho(q,p)$ is assumed to be smooth and positive it can be written as 
$$ 
\rho(q,p)= \exp \varphi(q,p)
$$
where the exponent must be such that $\int \rho(q,p) dq dp<\infty$ so that it  can be normalized to 1. Then from the steady-state condition
$$
\Frac{\partial \rho}{\partial t}=\Frac{\partial \rho}{\partial q}^{\top} \dot{q}+\Frac{\partial \rho}{\partial p}^{\top} \dot{p}=0
$$
 by  eliminating the (strictly positive) exponential on both members, it follows that
$$ \Frac{\partial \varphi}{\partial q}^{\top} p -\Frac{\partial \varphi}{\partial p}^{\top} V^2 q=\bmat \Frac{\partial \varphi}{\partial q}\\ \Frac{\partial \varphi}{\partial p}\emat^{\top} \bmat 0 &I\\-I &0\emat  \bmat\Frac{\partial H}{\partial q} \\ \Frac{\partial H}{\partial p}\emat=0
$$
where $H(q,p)$ is the Hamiltonian. The last equation is stating that the Poisson bracket   $\{ \varphi, H\}$ must be zero which means that $\varphi(q,p)$ must be a constant of the motion of the system. Now since $V^2$ is positive definite the only constant of the motion is, modulo a canonical  change of coordinates which does not interests us, a constant times the Hamiltonian. By integrability of $\rho$ we are lead to choose a negative  constant so that we can, with some hindsight, set $\varphi(q,p)= -\frac{1}{\beta} H(q,p)$ for some $\beta >0$. 
\end{proof}

\begin{rem} When $H(p,q)$ is a sum of two components $H(p)+H(q)$, the invariant distribution is {\em multiplicative}, i.e. $\rho(p,q)=\rho(p)\rho(q)$. 
\end{rem}
From now on, in order to avoid awkward notations such as $p(p)$, we shall change notation and denote the normalized momentum (i.e. the velocity random vector) by $\vb(\equiv p/m)$; where $\vb$ has sample space $ \Rbb^{3}$ and may change randomly from particle to particle. Hence:

\begin{cor}
The Maxwell-Boltzmann distribution $p(\v)$  is the {\em marginal} of the invariant distribution $\rho_{\beta} (p,q)$ with respect to the velocity components.
\end{cor}

Something which may look strange is that under the invariant distribution the velocity process turns out to be independent of the motion process $t\to \qb(t)$. This is however just a consequence of differentiability and can be proven in general for stationary processes using the spectral  representation.

 \nd Denote a nonsingular square root of $V^2$ by $V$. One can normalize and simplify the formulas of Hamiltonian mechanics   by {\em complexification} by defining\\
 $
 z:= p+j Vq,\,\text{where} \q j:=\sqrt{-1}, 
 $ 
then
$$
\dot{z}= -V^2q +jV\dot{q}= jV[ jVq +p]= jV z
$$
whereby  the phase space can be made into a complex Euclidean 3N dimensional space with $H(q,p)\equiv H(z)= \met \|z\|^2$. If $V$ is taken to be the unique symmetric square root of $V^2$  then the  Hamiltonian flow $\{\Phi (t)\}$ becomes  a continuous {\bf unitary} group  on the  $3N$-dimensional complex  Euclidean space. This will be generalized later on to more general phase spaces and become an helpful tool for the representation  of the Hamiltonian flow.
Under the invariant distribution  $z(t):= \Phi(t) z(0)$ becomes a complex   {\em stationary Gaussian process}.

\subsection{\bf Entropy and Probability}
There has been a century  long debate regarding the relations between Entropy and Probability, in particular between  Boltzmann famous formula  for Entropy $S=k\log W$ and probability. For the  large literature about this debate, see \cite {Jaynes-65} and  the list of  references in \cite{Davey-11}. We shall shortcut all of this debate and take a naive direct approach.

Basic observation: the integral $\Delta Q$  of the heat supply on a certain interval of time, or about a path describing a transformation in a thermodynamical state-space, must be equal to the increase of the {\em total internal (kynetic) energy} of the system (1st Principle). For an ideal gas the total energy of the gas is obtained by summing all individual kynetic energies of the particles.\\
The total internal energy $E(x)$ in the arrival macroscopic state $x$ of the transformation  is therefore the sum of the individual kynetic energies $\met m\v_i^2= \met m \sum_{k=1}^3 \v_{i,k}^2$ for $i=1,\ldots,N$. Assume that the particles move independently each with the   Maxwell-Boltzmann   probability distribution $p(\v)$ describing the random variables $\vb_k^2,\, k=1,2,3$, you get the identity  
$$
\Delta\,(\Frac{Q}{kT})=\sum_{i=1}^N \Frac{m  \v_{i}^2}{2kT}\,- E(x_0) = -\sum_{i=1}^N\log p(\v_i) +\text{const.}\simeq - N\int_{\Rbb^{3}} \log p(\v) \,p(\v)d\v + \text{const.}
$$
since  by independence  for $N$ very large the average term in the middle  can be approximated by expectation with respect to $p(\v)$, 
 \beq
 -N\,\E p(\vb)= -N\,\int_{\Rbb^3} \log p(\v) \,p(\v)d\v,
 \eeq
(the sample space of the random variables $\|\vb_{k}\|^2$ being actually  $\Rbb_{+}^{3}$). Hence the {\em specific entropy variation} (per particle)  is $-H(x)$ where
\beq \label{NegEntropy}
H(x)= k \,\int_{\Rbb^{3}} \log p(\v) p(\v)d\v\,.
\eeq
so $H=-S$ is {\em neg-Entropy} which should {\em decrease} for irreversible transformations.

\begin{rem} The neg-entropy \eqref{NegEntropy} depends on the macro state $x$ of the system since the Maxwell-Boltzmann distribution of a system in macro-state $x$ must describe the velocities of the particles which have a total energy $E(x)$ (microcanonical ensemble).  Hence the distribution is actually a {\bf conditional distribution, given $x$} which should (and will)  be denoted $p_{x}(\v)$.
\end{rem}
 Assume we are still describing a rarified  ideal gas. Then in any transformation the macro-state $x$   varies in time depending on the global motion  of the system, which is described probabilistically  by a trajectory $t\to \vb(t)$ of the {\bf  stationary stochastic process} $\vb(t)$. In fact, for a rarified  ideal gas, $x(t)$ only depends on the squared norm of the random vector $\vb(t)$.\\
  Suppose a transformation drives the system from state $x_0$ to some other state $x$ at time $t$, then 
$$
H(x)-H(x_0)=k \int_{\Rbb^3} [\log p_{x}(\v)- \log p_{x_0}(\v)]\,p_{x_0}(\v)d\v\,= k \,\E  \log \Frac{p_{x}(\vb)}{p_{x_0}(\vb)}
$$

\begin{rem}[Historical remark]
 Boltzmann apparently discovered a  connection between his (non-probabilistic) formula $S=k\log W$ valid for systems with a discrete configuration space, with the formula \eqref{NegEntropy}; he used a $\rho\log\rho$ formula as early as 1866  interpreting $\rho$ as a density in phase space—without mentioning probability—but since $\rho$  satisfies the axiomatic definition of a probability density we can retrospectively interpret it as a probability anyway. J. W. Gibbs gave an explicitly probabilistic interpretation in 1878. 
\end{rem}
In the following lemma  we point out the self-evident fact that, for a finite ensemble of particles, the  velocity process must be  a rather trivial kind of stochastic process. It must in fact be  a {\em purely deterministic (or non-regular) process}. This notion  goes back to Wold, Kolmogorov and Cram\`er. See \cite{LPBook} for an up to date review of the concept.
\begin{lem}
 For an isolated ideal gas system of $N$ particles in thermal equilibrium, the microscopic state process
$$
\xib(t)= \bmat \qb_1(t) &\qb_2(t) &\qb_3(t) &\vb_1(t) &\vb_2(t) &\vb_3(t) \emat^{\top}\,\in \Rbb^{6}
$$
 is a {\bf purely deterministic} stationary Markov process
and  the velocity vector  $\vb(t)$ is a memoryless function of the state; i.e.  $\vb(t)= \bmat 0_3& I_3\emat \xib(t)$. Hence the velocity process is also purely deterministic.
\end{lem}

\subsection{\bf On the Proof of Entropy Increase}
Consider  now  an   ideal gas system of $N$ particles performing some  transformation  driving the system  from macro-state $x_0$ at time $t=0$ to some other macro-state $x$ at time $t$.
For any such  general transformation the entropy should  increase. How do we prove this?

We want a proof based on the  ``microscopic" model of the system. That is,  based on the transition $\vb(t_0)\to \vb(t)$ which causes  the transition  $p_{x_0}(v)\to p_{x}(v(t))$.

Because of the Markovian (although  purely deterministic) character of the velocity  we can imagine $\vb(t)$ for $t\geq 0$, still being described by the initial probability $p_{x_0}(v)$.

Assume the relation $x\to p_{x}(\cdot)$ is 1:1 (in Statistics: the model is identifiable), in fact we only need this property locally  about $x_0$; then 
$$
\int_{\Rbb^{3}} \,\log \Frac{p_{x_0}(v)}{p_{x }(v)}p_{x_0}(v) \,dv= \E_{x_0} \log \Frac{p_{x_0}(\vb)}{p_{x }(\vb)}=K(x_0,\,x)
$$
is the  {\em Kullbak-Leibler (pseudo)-distance}   which is always {\em nonnegative} and can be zero only if $x=x_0$.\\
The neg-entropy variation of the {\bf random} particle going from $x_0$ to $x$ along a chain of macroscopic states $\{x_0\to x_1 \to \ldots \to x_n\equiv x\}$ is the expectation of the random difference (note the bold characters)
$$
\Hb(x)-\Hb(x_0)= k\,  \log \Frac{ p_{x}(\vb )}{p_{x_0 }(\vb_0)}= k\sum_{j=1}^n\,( \log p_{x_j}(\vb_j)-\log p_{x_{j-1} }(\vb_{j-1})\,)
 $$
{\em If $\vb$ was an ergodic process}, then  by the ergodic theorem, 
 $$ 
\lim_{n\to \infty} \,\Frac{1}{n}\;  \sum_{j=1}^n\,  \log \Frac{ p_{x_j }(\vb_j)}{p_{x_{j-1} }(\vb_{j-1})}= \E_{x_0}\log \Frac{ p_{x }(\vb)}{p_{x_0 }(\vb)}= -K(x_0,\,x)<0\,,
$$
for every $x_0$, with probability 1. So if the $\vb$ process was ergodic, the neg-entropy would decrease! Obviously (and unfortunately) Boltzmann could not know the ergodic theory of stochastic processes which was developed much later by von Neumann, Birkhoff and Kolmogorov in the 1930's and could only use intuitive arguments. In fact, since it is purely deterministic, {\em the velocity process  $\vb(t)$ of a finite  ensemble of  gas molecules  cannot be ergodic!}

\nd Ernst Zermelo \cite{Zermelo-896}, \cite{Zermelo-896b} pressingly criticized Boltzmann argument by an irrefutable   reference to {\em Poincar\`e recurrence theorem}: 
\begin{center} 
{\em Any function of the state of an  isolated system of point masses constrained in a finite volume  must cycle (have  a periodic behavior).} 
\end{center}
Hence the velocity process $\vb(t)$ must be periodic (which checks with its purely deterministic character). So, again,
{\bf Periodic processes cannot be ergodic!}  
hence $\vb(t)$ cannot be ergodic! The previous limit argument does not work.
\medskip

Boltzmann worked very hard to modify his neg-entropy model by including collisions among particles. But his famous  equation and the   proof of the related {\it H-Theorem} are not based on a microscopic mechanical  model and Kynetic theory. They just describe density functions. From this point of view the argument  seems still to be unsatisfactory. See e.g \cite[p.  164]{Prigogine-80}.

\section{\bf   Coupling to an  infinite dimensional heat bath} 

 Since the early attempts in the literature, the transition from the kynetic to  macroscopic thermodynamical models of nature has been studied and developed inspired by the idea that one is really considering  a  {\em limit situation} when the number of particles of the underlying system is ``very large" (say compared to  Avogadro's number) and in the background one can imagine a passage to    a ``thermodynamic limit" when $N\to \infty$, that is, the number of degrees of freedom of the system goes to infinity. The mathematical details of this passage to the limit are often left obscure however. One is lead to think that we should perhaps start by considering systems with $N=\infty$ from the outset.
Infinite-dimensional descriptions are however often criticized a being ``non physical". In this context one may refer to an often quoted sentence of  Jan C. Willems \cite{Willems-86}: {\it The study of infinite-dimensional [say ``non-complete"  in his language] systems does not fall within the competence of system theorists and could be better left to cosmologists and theologians..}

In this paper we shall nevertheless depart from this suggestion and attempt to show that a first-principle microscopic  understanding and modeling of  irreversibility, which  macroscopically is empirically described as  friction, resistance, dissipation etc.  may be achieved by going to   infinite dimensions. What we have just learned in the previous section is that the microscopic evolution of a finite number of mechanical particles or systems cannot show a truly stochastic (and therefore irreversible) character. It seems indeed  that these phenomena cannot be explained from first principles by finite-dimensional microscopic dynamical  models alone.   We shall instead attempt to   explain, at least in a linear framework,   the origin and the rationale of stochastic modeling of physical systems, and consequently of macroscopic irreversibility, by assuming a coupling  of the (necessarily  lossless)   deterministic  model of a   physical system to an  {\em infinite-dimensional environment}, which is called by physicists a {\em heat bath}. In the next paragraph  shall very quickly make reference to a related mathematical field which   has been object of an intense study in the past decades and may be related to our program. It may in fact   also involve infinite dimensional dynamics.

\subsection{Linear Chaos }
Tthere are    concepts like  ``Chaos" and "chaotic dynamics" which have been introduced  to study certain deterministic dynamical systems which can show a stochastic-like behavior. See    \cite{Devaney-89} for  a readable review.  In the literature  chaotic dynamical systems are however almost invariably   assumed to be {\bf non-linear} finite-dimensional  dynamical systems evolving on a compact phase space endowed with a  hyperbolic structure.   In all  such cases, the interesting invariant sets on which irregular "stochastic" behavior occurs, are "small", typically of zero Lebesgue measure.The outputs of such systems, on such invariant sets are effectively  "random"  processes but only of the {\em finite--state}
type, e.g. Bernoulli i.i.d. processes. One would instead like to understand how to generate stationary processes with a continuous state space, say Gaussian processes taking values in
$\Rbb$ or $\Rbb^n$.\\
Curiously enough, the possibility  of chaos in {\em linear} systems and its
nature has not been reputed an interesting object of study  in the mathematical
community until  recently.  There are now some  characterizations of
  chaos  (see \cite{Herrero-92}, \cite{Godefroy-S-91}) which make it possible, specifically for infinite-dimensional linear systems. Based on this research  it can safely be stated that   linear chaos can happen only in infinite dimensions \cite{Grosse-M-011}. These conditions support the 
conclusion that, at least in a linear setting,  a purely non-deterministic (ergodic) stochastic behavior can only  come about in an infinite-dimensional phase space.  \\
  Since we want to understand irreversibility of the motion of finite-dimensional ``theoretically conservative" system (say a  Newtonian object moving in open space). The key structural fact   that we need to add to the picture is  the {\em coupling} of  the system  with an {\em infinite-dimensional} surrounding environment. This environment  we call  the {\em heat bath} and is the designated cause  which generates drag, resistance or energy loss. Physically, the surrounding environment is always present. The analysis of this coupling is however highly non trivial and in the   literature there have ben several attempts  to understand it  from different perspectives  \cite {Ford-K-M-65} \cite{Lewis-T-75} some of which not very convincing from our point of view.\\
Understanding the coupling  mechanism   subsumes understanding  the departure from the conservative lossless behavior of the system  and eventually how it may become (or look like)    a dissipative, stochastic one. For a long time this has been a general goal in some physical literature. \\
In the approach of this paper  which is believed to be original,  we fill in  a new unexpected ingredient, namely {\em  feedback control theory}. Using  just some basic ideas will, in our opinion, clarify in a neat way the mechanism of this transition. \\
Three typical  examples of analysis of this coupling will be presented below.

\subsection{\bf Example 1: The Infinite Lossless Electrical  line}



\begin{center}
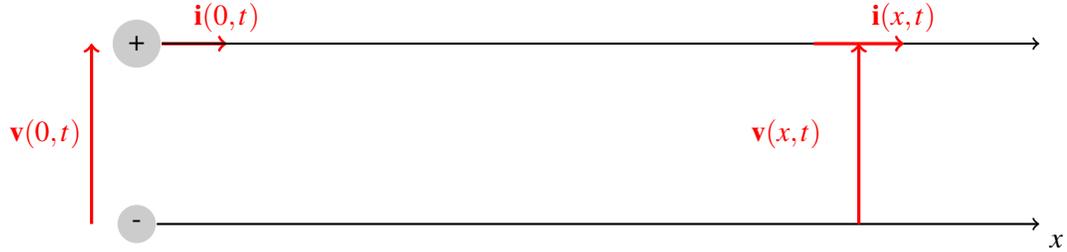
\begin{figure}[h]
\begin{tikzpicture}[scale=1.2]

\draw[red, very thick,->] (-2,0) -- (-1,0)node[above] {$\ib(0,t)$ };
\node at(-2.5,-1)[left, red] {$\vb(0,t)$};

\draw[thick,->] (-2,0) -- (8,0);
\draw[red,very thick, ->] (5.5,0)--(6.5,0) node [above] {$\ib(x,t)$ };
\draw[thick,->] (-2,-2) -- (8,-2)node[below right] {$x$ };

\node at (-2,0) [circle,draw=black!1,fill=black!20]{+};
\node at (-2,-2) [circle,draw=black!1,fill=black!20]{-};

\draw[red,very thick,->] (-2.5,-2) -- (-2.5,0);

\draw[red,very thick,->] (6,-2) -- (6,0);
\node at(5.7,-1)[left, red] {$\vb(x,t)$};
\end{tikzpicture}
\label{figure1}
\caption{A semi-infinite electrical line}

\end{figure}
\end{center}

\vspace{-0.3 cm}

Let $L$ and $C$ be inductance and capacitance of the line per unit length, then the voltage and current on the line at distance $x$ from the origin,  $\v(x,t)$ and $i(x,t)$, satisfy 
 \beq \label{WaveEl}
\bmat
\frac{\partial v}{\partial t}\\ \frac{\partial i}{\partial t}\emat  =  
 \bmat 0 & \frac{1}{C}\,\frac{\partial}{\partial x}\\    \frac{1}{L}\,\frac{\partial}{\partial x}& 0 \emat    \,  \bmat v \\ i\emat    \q \equiv \q
 \begin{cases}
 \frac {\partial \vb}{\partial t} &= \q \dot{\vb} \\
  \frac {\partial \dot{\vb}}{\partial t} &=  \frac{1}{LC}\,\frac {\partial^2 \vb}{\partial x^2}  \end{cases} 
\eeq
  the second of which is   in Hamiltonian form. They  both  lead to the classical wave equation
\footnote{N.B.  Note that the second derivative operator in $L^2$ is (symmetric) negative semidefinite. Therefore the analog of the finite dimensional  potential matrix $V^2$ is the operator $- \frac {\partial^2}{\partial x^2}$.
}
\beq \label{HamiltWave}
  \frac {\partial^2 v(x,t)}{\partial t^2} - \Frac{1}{LC}\,\frac {\partial^2 v(x,t)}{\partial x^2}=0,\q  \frac {\partial^2 i(x,t)}{\partial t^2} - \Frac{1}{LC}\,\frac {\partial^2 i(x,t)}{\partial x^2}=0   \,.
\eeq
  It is well-known to electrical engineers that the input impedance seen from the terminals at $x=0$ of such  an  Infinite lossless line is {\bf purely resistive} with impedence
$$
Z_0:=\frac {v(0,t)}{i(0,t)} = \sqrt  { \frac{L}{C}}\q\q \text{Ohms}\,.
$$
Hence from the input terminals the infinite lossless line behaves like a {\bf dissipative system}. This looks quite surprising. How could this be? There are no resistances so  the system is conservative. We make no assumptions on how the infinite line is terminated (no radiation due to coupling with an external world) so no dissipation is visible. In fact if we close the line at any finite length the behavior becomes {\em purely oscillatory!}.

This is actually a manifestation of a general fact. To  understand the phenomenon in general, consider   any lumped {\em conservative} $L C$
circuit,  composed by an arbitrary (but finite) number of linear 
capacitors and inductances connected at the left-end point, situated at $x=0$, to a semi--infinite electrical lossless  line  like that in fig \ref{figure1}.
The overall system is still infinite-dimensional and conservative. However after some  analysis we shall discover the following
\nd{\em Fact:} If we connect to the lossless lumped circuit (still connected to the  the infinite  line)  a linear observer, we end up observing a signal which is described by a linear finite-dimensonal equation which has a {\em dissipative dynamics} and a   ``white noise-like" input (although deterministic). Further,  the model has an asymmetric {\em irreversible}  behavior, that is: {\em The model changes if we change the direction of time~!!}

\subsection{\bf Some Analysis}

First we need to identify the boundary conditions:
At $x=0$ the current $i_{o}(t) := i(0,t)$ acts on the electrical load
inducing a constrained motion  with voltage $v_{o}(t) :=
v(0,t)$, which must be related to the acting current by the  impedance of the load.
The whole system is assumed to be in steady state and we shall consider its evolution for $t\in \Rbb$.  Using  (formal) Fourier transforms denoted by hatted symbols,   one can write
  \begin{equation}  
  \hat{v}_{o}(j\omega) = Z_{o}(j\omega) \hat{i}_{o}(j\omega)     \label{IMPEDel}  
\end{equation}  
 where $Z_{o}(j\omega)$ is the  impedance of the load seen from the connecting
point with the line at $x=0$. Equivalently, the variable $j\omega$ could be substituted by a bilateral  Laplace transform complex variable $s$.

\nd Note: Here to justify the use of the standard $L^2$ Fourier transforms one should  assume that $\v_0(t)$ and $i_0(t)$ are (finite energy) signals in $L^2(\Rbb)$. At this point this is not obvious and will come out as a byproduct of Theorem \ref{MainElThm}  later on; for now  one should start by using  generalized functions and generalized Fourier transform, see e.g. \cite{Lewis-M-84}.
However since the mathematics would become more involved  and obscure the argument, we shall not adventure into it and proceed formally.

The composite system can be described by the following state equations  
\begin{align}
\bmat
\frac{\partial v}{\partial t}\\ \frac{\partial i}{\partial t}\emat & =  
 \bmat 0 & \frac{1}{C}\,\frac{\partial}{\partial x}\\    \frac{1}{L}\,\frac{\partial}{\partial x}& 0 \emat    \,  \bmat v \\ i\emat             \label{WAVE2}  \\     
\dot{\xi}(t) & =   A\xi(t) + b_{o}i_{o}(t)                \label{STATEel}  \\
  \v_{o}(t)& =  c_{o}\xi (t)                     \label{VOel}  
  \end{align} 
where $\xi(t)$ is a $n$ dimensional state variable of the load, where $i_{o}(t) \equiv  i(0,t)$ acts as in {\em input function} and similarly $\v_{o}(t) \equiv \v(0,t)$ as an {\em output}. To adhere to standard notations in System Theory, later  we shall name  them $u(t)$ and $y(t)$ respectively.\\
 The last two equations (\ref{STATEel}) and (\ref{VOel})   will then constitute a state-space  realization of the transfer function 
 \begin{equation} 
 Z_{o}(j\omega): = c_{o}(j\omega I - A)^{-1}b_{o}         \label{ZOel}   
 \end{equation}
which is in fact just  the electrical impedance $Z_{o}(j\omega)$ of the load. The realization will be assumed to be   a (minimal) realization in the sense of System Theory. \\
N.B. (important observation)  Since by assumption $Z_{o}(j\omega)$ is a {\em Lossless} impedance function, there
cannot be a direct feedtrough term in the realization \cite{Anderson-V-73}.

  We now assume suitable units have been chosen to insure $LC = 1$ (so
that the speed of propagation along the line is one). The evolution of the
composite system (3-4-5) can  then be seen as the evolution of a
conservative Hamiltonian system $\dot{z} = Fz $ with overall state (phase) vector
\begin{equation}
 z:= \left[\begin{array}{c} \xi \\ \v \\i \end{array} \right] \label{Z} 
\end{equation}  taking place in the phase space ${\bf H} := {\Rbb}^{n}\oplus
L^2_{2}({\Rbb_+})$. This vectorspace can be given a Hilbert space structure by
introducing the {\em energy norm}  
\begin{equation}
 \| \bmat \xi \\ v \\i \emat \|^2 = 1/2 \xi^{\top}\Omega
\xi + 1/2 \int^{+\infty}_{0} (C v^2 + L i^2) dx   \label{NORM}  
\end{equation}
where $\Omega$ is a symmetric nonnegative matrix representing the total energy
(hamiltonian) of the load, a quadratic form in the state $\xi$. By choosing $\xi$
minimally we can always guarantee  $\Omega >0$.

\begin{rem}[On finite total energy]
{\em
The above  Hilbert space  structure implies that the time  evolution of the system, in particular  of the infinite line, should occur with finite total energy. To stay within the $L^2(\Rbb)$ structure in space we just need to assume that at each fixed time $t$,  the energy of the line, at distance $x$, i.e. $1/2   (C \v(x, t)^2 + L i(x,t)^2) dx $ is  decreasing with the distance $x$ fast enough to make the integral finite. One may question if this is physically reasonable; in particular could    an infinite line have finite total energy.

Although mathematically consistent this assumption may seem to be physically unreasonable. 
 A physically sounder framework, could be to assume the time evolution to take place on a  space of functions (or distributions) which are just   locally $L^2$, e.g. the dual of some Sobolev space. This may be done in such a way as to keep formal similarity with the essentials of the $L^2$ mathematical structure but  would, on the other hand, make the whole discussion obscured by technical details so we shall not dwell on  this point further.}
\end{rem}

It is not hard to see that the $F$ operator  determining   the dynamical equations of the overall conservative system (3-4-5) must be  skew-adjoint on its natural domain (of
square summable functions satisfying the boundary conditions (\ref{IMPED})) and generate
an  energy preserving (i.e. orthogonal) group on the overall phase space ${\bf H}$. In particular in steady state the voltage  and current variables  $\v(\cdot,t)$ and $i(\cdot,t)$ satisfying the wave equation \eqref{HamiltWave}  should be absolutely continuous functions of the time variable $t$ evolving in the joint phase space  $\Rbb^{n}\times L^2(\Rbb_+)$.

\subsection{\bf The wave picture}
By normalization, both the voltage
and  the current of the line subsystem  obey  the same wave equation and can both   be expressed by the well known d'Alembert formula
\begin{equation}
 \varphi(x,t) =a(x+t) + b(x-t) \hspace{1.cm} t \in {\bf R},\q x \geq 0 \label{AB}
\end{equation}
 whereby the the group $\Phi(t)$, solution of the Hamiltonian equation \eqref{HamiltWave} is seen to be a combination of the forward and backward translation operator acting on the initial  conditions of the line at  time  $t=0$. The functions $a(x+t)$ and $b(x-t)$ are  called the {\em incoming} and {\em
outgoing} waves respectively.    \\
In order to get   symmetric formulas introduce the {\em charge function} of the line at time~$t$:
$$
q(x,t)= \text{electrical charge of the section}\; [0,\;x]
$$
Then since
$$
i(x,t)= \frac{d q(x,t)}{dt}\,;\qquad \v(x,t)dx=Cdq(x,t) \Rightarrow \;  \v(x,t)= C \,\frac{d q(x,t)}{dx}
$$
we obtain (for $C=1$)
\beq  \label{vfel}
\v(x,t)= \frac{d \,q(x,t)}{dx}\,=a'(x+t) +b'(x-t)\,; \quad  i(x,t)= \frac{d q(x,t)}{dt}\,= a'(x+t)- b'(x-t)   
\eeq
where the primes stand for derivatives.

Now,  following  a standard idea of Scattering Theory, as used e.g. in \cite[p. 114]{Lewis-T-75}, one can  also solve \eqref{vfel}, for the functions $a'$ and $b'$  obtaining
\begin{subequations} \label{ICe}
\begin{eqnarray}   
2 a^{\prime}(x+t)&=&\v(x,t) + i(x,t)  \\  
2b^{\prime}(x-t) &=& \v(x,t)-i(x,t)           
\end{eqnarray}  
\end{subequations}
so that, in the present setup the functions $a'$ and $b'$ which  enter in the scattering representation of the
state vector $(\v,i)^{\top}$ are determined  by the initial data  at   time zero , that is by $\v_{0}(x) := \v(x,0)$ and $ i_{0}(x) := i(x,0)$. In fact, putting $t=0$ in \eqref{ICe}, one obtains 
\begin{eqnarray}  
2 a^{\prime}(x)&=&\v_0(x) + i_0(x)  \label{ICa}  \\ 
2b^{\prime}(x ) &=& \v_0(x)-i_0(x) \,. \label{ICb}
\end{eqnarray}  
This  system of equations determines 
$a'(x)$ and $b'(x)$  on the half line  $x \geq 0$  as  the initial conditions $v_0(x)$ and $i_0(x)$ are only defined as elements of $L^2(\Rbb_{+})$ while on the complementary half lines the functions are not defined \footnote{Note  however that the support of the functions depends only on the choice of the initial coordinate at $x=0$ of the line which is inessential and could in fact be moved arbitrarily in space.}. It should also be remarked   that,  since  the  initial data functions $[\v_0(\cdot),i_0(\cdot)]^{\top}$ may be  arbitrary functions in $L^2_{2}(\Rbb_{+})$ the  derivatives $a'|_{x\geq 0}$ and $b'|_{x \geq 0}$ are also arbitrary and may well be called   "free variables" (i.e. inputs) in the sense of J.C. Willems \cite{Willems-86}.

The first members of  \eqref{ICa} \eqref{ICb}  propagate  in   time,  the first inheriting from $a(x+t)$ the character of an {\em incoming wave} while $b^{\prime}(x,t)$ having dually the character of an {\em outgoing (or reflected) wave}.
To make this aspect more evident we shall introduce some notations. \\
Temporarily denote by boldface symbols the elements of the Hilbert space $L^2(\Rbb_{+})$ so that the initial data $v_0(x)$ and $i_0(x)$ as functions  in $L^2(\Rbb_{+})$ are denoted $\vb_0$ and $\ib_0$. Similarly, write $\ab_0^{\prime}:=\met (\vb_0 +\ib_0)$ and $\bb_0^{\prime}:=\met (\vb_0 -\ib_0)$. Now  introduce the {\em translation group} $\{\Sigma(t); t\in \Rbb\}$ on $L^2(\Rbb)$ which acts by shifting the argument of a function by $t$:
\beq
[\Sigma(t)f](x):= f(x+t),\qq[\Sigma^{*}(t)f](x):= f(x-t)\qq f \in L^2(\Rbb)\,. \notag
\eeq
Then we can represent the above outgoing and incoming waves as a result of {\em translation in time} of the initial data, i.e.
\beq \label{TranslICa}
a^{\prime}(x+t)= [\Sigma(t)\ab_0^{\prime}](x),
\eeq
which could be compactly rewritten  as $\ab_{t}^{\prime} = \Sigma(t)\ab_0^{\prime}$.

Since we shall need information about the behavior of the $b^{\prime }$ function for negative arguments we shall use a dual relation to  \eqref{TranslICa}, namely
\beq \label{TranslICb}
 b^{\prime}(-x-t)= [\Sigma(-t)\bar\bb_0^{\prime}](x),\qq \text{\rm where}\qq \bar\bb_0^{\prime}(x) = \bb_0^{\prime}(-x)
\eeq
which can be written compactly as $\bb_{-t}^{\prime}= [\Sigma(-t)\bar\bb_0^{\prime}]$.\\
Equations \eqref{TranslICa} and \eqref{TranslICb}  make explicit the ``movement" of the initial conditions by translation in time. Note that the support of $\ab_{t}^{\prime}$ is shifted from $\Rbb_+$ to the half line $[t , \,+\infty)$ while that of $\bb_{-t}^{\prime}$ is changed from $\Rbb_-$ to the half line $(-\infty, \,t]$. 

 Consider now the two waves observed at the endpoint $x=0$, denoted  $ a^{\prime}(t)$ and $b^{\prime}(t)$. These are signals which live for all $t\in \Rbb$ since the right members of \eqref{ICe} are steady-state solutions of the wave equation  which describes the system for all times (after the two waves have travelled for an infinite amount of time).  They evolve in time  by  the action of the left and right time translation operators acting on the initial condition $\ab_0' \in L^2(\Rbb_{+})$ and $\bar\bb_0' \in L^2(\Rbb_{-})$    which are determined by the known initial values of current and voltage of the semi-infinite line. We shall show that these two signals are related by a {\em scattering function}, a unitary transfer function which in our case is in fact a rational all-pass function.

Next by putting $x=0$ in \eqref{ICe}  one obtains  the identities  for time-dependent variables
\begin{eqnarray}
\v_0(t)\equiv \v(0,t) & = & a'(t) + b'(-t) \label{FSCATTe}\\    
i_0(t)\equiv i(0,t) & = & a'(t) - b'(-t)   \label{SCATTe}  \qq t\in \Rbb
\end{eqnarray} 
 which are  the "steady--state" boundary condition at $x=0$, relating $\v(t)$ and
$i(t)$ constrained by the realization of the impedance (\ref{IMPEDel}). It is immediate to check that
(\ref{FSCATTe}) and (\ref{SCATTe}) are both interpretable as {\em state  feedback}
laws on the conservative load system $(A,b_{o},c_{o})$. In fact, by identifying the current $i_0(t)$ as the system input,  they  can be written   
\begin{eqnarray}
 i_0(t) &=& -\v_0(t) +2a'(t) = - c_0\xi(t)+2a'(t)  \label{FFel}\\       
i_0(t) &=& \v_0(t) -2b'(-t)    = c_0\xi(t)    -2b'(-t)   \label{BFel} 
\end{eqnarray}
which  appear respectively as a {\bf negative} and {\bf positive state feedback} relations  applied to the load system \eqref{STATEel}, \eqref{VOel}. A further change of notations:
$$
 w(t) : \equiv a^{\prime}(t),\qq  \bar w(t) : \equiv -b^{\prime}(-t)
 $$
 helps to streamline the formulas and eventually will make contact with their stochastic version to be made explicit in part two of this paper.  With them, the feedback relations \eqref{FFel} and  \eqref{BFel} give rise to the following pair of representations of the state dynamics of the load :
\begin{eqnarray}
 \dot{\xi}(t) & = & (A -b_{o}c_{o})\xi(t) + 2b_{o}w(t) \label{FRel}\\   
 \dot{\xi}(t) & = & (A + b_{o}c_{o})x(t) +2b_{o}\bar w(t)  \label{BRel} 
\end{eqnarray}
which should be coupled to the output equation $\v_0(t)= c_0\xi(t)$ of the load. Introducing the feedback matrices 
\beq \label{Gammas}
\Gamma:= (A -b_{o}c_{o})\, \qq \bar \Gamma:= (A + b_{o}c_{o})
\eeq
we obtain two dual expressions for the transfer functions from the incoming (resp. outgoing) wave to the output $\v_0(t)$ of the load
\beq
\hat y(s)= 2 c_0(sI-\Gamma)^{-1}b_0\,\hat w(s)= 2 c_0(sI-\bar\Gamma)^{-1}b_0 \,\hat{\bar{w}}(s)
\eeq
where the hat stands for (doble-sided) Laplace transform\footnote{ Equivalent to Fourier transform.}. After eliminating $\hat y(s)$ we get
\beq \label{Kfrac}
   \,\hat{\bar{w}}(s)= \Frac{c_0(sI-\Gamma)^{-1}b_0} {c_0(sI- \bar{\Gamma})^{-1}b_0} \;\hat w(s):=  K(s)\hat{w}(s)
\eeq
which turns   out to be  the scattering function relating the two waves.

A well known dictum of contro theory states that  feedback does not modify the zeros of the transfer function $Z_0(s)$. It follows that the
zeros of the rational functions $c_0(sI-\bar\Gamma)^{-1}b_0$ and   $c_0(sI-\Gamma)^{-1}b_0$ coincide and  cancel out in forming the quotient \eqref{Kfrac}. Hence the poles of $K(s)$ turn out to be the eigenvalues of $\Gamma$ while the zeros turn out to be  the eigenvalues of $\bar \Gamma$. We shall see below that these eigenvalues are in fact  opposite  of each other.

\begin{thm} \label{MainElThm}
Assume the realization (\ref{ZOel}) is minimal. Then  the dynamical sytems (\ref{FRel}) and (\ref{BRel}), respectively,
are asymptotically stable and antistable, in fact,
  \begin{equation}  
\Re e\,\lambda(A - b_{o}c_{o}) <0,\hspace{1cm}  \sigma(A + b_{o}c_{o}) = - \sigma(A -
b_{o}c_{o})   \label{EIGel} 
 \end{equation}  the symbol $\sigma(A)$ denoting the spectrum of the matrix $A$. \\
Moreover the two representations are related by a kind of  "change of white noise input"
formula, of the type \eqref{Kfrac}  where  $K(s)$ is precisely the scattering function associated to the boundary
condition (\ref{IMPEDel}), i.e.   
\begin{equation} \label{KfromZ} 
 K(s) = \frac{Z_{o}(s) -1}{Z_{o}(s) + 1}\,.  
 \end{equation} 
 In fact, $K(s)$ is a rational  {\bf inner function}. 
 \end{thm}

\begin{proof}
The representation \eqref{KfromZ} follows from \eqref{FFel} and \eqref{BFel} which can be rewritten in terms of Laplace transforms as $Z_0(s)\hat {i}_0(s)+\hat {i}_0(s)= 2 \hat w(s)$ and as  $Z_0(s)\hat {i}_0(s)-\hat {i}_0(s)= 2 \hat {\bar{w}}(s)$, leading to   
\begin{eqnarray}
\hat{\v}_0(s) & = & (1+ Z_{o}(s))^{-1}2 \hat w(s)   \\
\hat{\v}_0(s) & = & (1- Z_{o}(s))^{-1} 2 \hat{\bar{w}}(s)  \label{TF}
\end{eqnarray}
which  implies \eqref{KfromZ}. Since $Z_{o}$ is {\em Positive Real} the
same is true for $(1+ Z_{o}(s))^{-1}$ which must therefore have all poles
strictly inside the left half plane. By minimality this implies that the
first condition in (\ref{EIGel}) is in fact true.

 Expressing the impedance of
the load  as $ Z_{o}(s)= N_{o}(s)/D_{o}(s)$ we conclude from the
argument above that all zeros of the polynomial $D_{o}(s) + N_{o}(s)$
must lie strictly inside the left half plane. Now it is well known that a rational 
lossless impedance must be the ratio of even and odd polynomials in $s$ \cite{Foster-24}.
From this it is not hard to see that the zeros of the denominator 
$D_{o}(s) - N_{o}(s)$, of the transfer function $ (1- Z_{o}(s))^{-1}$ are
 just the opposite of the zeros of $D_{o}(s) + N_{o}(s)$.
This proves the second relation in (\ref{EIGel}).

 The last statement of the
theorem follows now easily from the above, after eliminating $\hat {\v}_0$ in
(\ref{TF}).  
\end{proof}

 Next we want to study the time evolution of  an {\bf arbitrary linear observable}, whose  value say $y(t)$  may   well be called an {\em output variable} of the load system. Any such observation channel being  described by an arbitrary linear combination of the state of the load $\xi$ with perhaps a direct
feedthrough term from the current  of the line $i_{o}$,  say:
\begin{equation}
 y(t)= c\xi(t) + di_{o} (t)     \label{OUTPUTel}  
 \end{equation}
where $c$ is an $n$-dimensional row vector. This is a quite general  linear functional of the state of the  conservative circuit coupled to the infinite-length electric line. We shall show that, stated in thermodynamical language, the infinite line  will  be  playing the role of    a {\em heat  bath} coupled to the finite-dimensional lossless circuit and that 
this coupling will induce the famous double action of both dissipation and stochastic behavior  of the output function.  This is described in physics as  a (generalized)  {\em Langevin equation}. In fact it will be a generalized Langevin-like vector equation of dimension $n$.

A combination of the state equations \eqref{FRel}, (\ref{FFel}) with \eqref{OUTPUTel}    yields the
following  representation of the output signal $y$ :
\begin{subequations}\label{Langevin}
\begin{eqnarray}
 \dot{\xi}(t) & = & (A -b_{o}c_{o})\xi(t) + 2b_{o}w(t) \\ 
   y(t) & = & [c- dc_{o}]\xi(t) + 2w(t)   \label{FRout}
\end{eqnarray}
\end{subequations}
 while, dually, a combination of the state equation \eqref{BRel}, (\ref{BFel}) with \eqref{OUTPUTel}  yields
 \begin{subequations}\label{LangevinBar}
\begin{eqnarray}
 \dot{\xi}(t) & = & (A + b_{o}c_{o})\xi(t) +2b_{o}\bar{w}(t)  \\
 y(t) & = & [c+ dc_{o}]\xi(t) + 2\bar{w}(t)  
\end{eqnarray}
\end{subequations}
which are two dynamical descriptions  of the (same !) output signal $y$  by two linear {\em finite dimensional } state equations   driven either by the outgoing scattering wave $ a^{\prime}(t)$ or by  the incoming $b^{\prime}(-t)$
which may again be looked upon as {\em free inputs} as they are determined by arbitrary initial distributions of voltage and current along the line. We shall actually take the liberty of calling them {\em noise variables} and then, with some foresight, call \eqref{Langevin} and \eqref{LangevinBar}  {\em forward and backward generalized deterministic Langevin equations}.

 In Part Two of this paper we shall show  that these two  linear models  have in fact the same  structural
properties of a  {\em forward-backward pair of   stochastic models} driven by a  {\em forward-backward pair of white noise processes} exactly of the kind first introduced in \cite{Lindquist-P-79} and extensively  studied in \cite{LPBook}. Note incidentally that {\em the inner scattering function is independent of the particular observable}, i.e. of  the specific output equation.  In fact,
Introducing the notations 
$$
h:=c- dc_{o}\,\q \bar h:= c+ dc_{o}
$$
we can write the outputs of \eqref{Langevin} and \eqref{LangevinBar}  as $y(t)= \int W(\tau) w(t-\tau) d\tau$ and   $y(t)= \int \bar W(\tau) \bar{w}(t-\tau)d\tau$ where the transfer functions
\begin{align}
\hat W(s)&= 2\left[ h(s I-\Gamma)^{-1})b_0 +1\right]\\
\hat {\bar W}(s)& = 2\left[ \bar h(s I-\bar \Gamma)^{-1})b_0 +1\right]
\end{align}
 assuming minimality, are transfer functions having all   poles at the eigenvalues of $\Gamma$ and $\bar \Gamma$ respectively. Again, by the invariance of the zeros under linear feedback, it can be  checked by standard state-space calculations that 
\beq \label{InnerW}
\hat {\bar W}(s)^{-1}\hat W(s) = K(s)
\eeq
where the hats stand for (double sided) Laplace transforms and the {\em  inner function $K$  is the same as the original scattering function of Theorem} \ref{MainElThm}. In other words the scattering function is invariant with respect to the choice of the linear observable $y$.
In fact it is easy to check that  the above calculation applies verbatim to the evolution of {\em any} linear  observable, (or output) $y$ of the system, described by a linear functional of the type (\ref{OUTPUTel}).

\begin{rem} This  dual forward-backward representation phenomenon is very similar to what happens    to linear  stochastic realizations of a stationary process.   One should note that  the evolution of
$y(t)$ {\em backwards in time} is governed by a dynamical  model
(\ref{BRel}) which does {\em not} correspond to the trivial {\em reflexion} of
time transformation $t \rightarrow -t$ on (\ref{FRel}). The latter  is in fact not quite a {\em change of direction}. Hence the difference of temporal evolution of the variable $y(t)$ in the reverse direction of time, resembles a truly
stochastic phenomenon which has no counterpart in finite dimensional
deterministic systems.
\end{rem}

\subsection{\bf Example 2: The Mechanical string}\label{MechSubs}

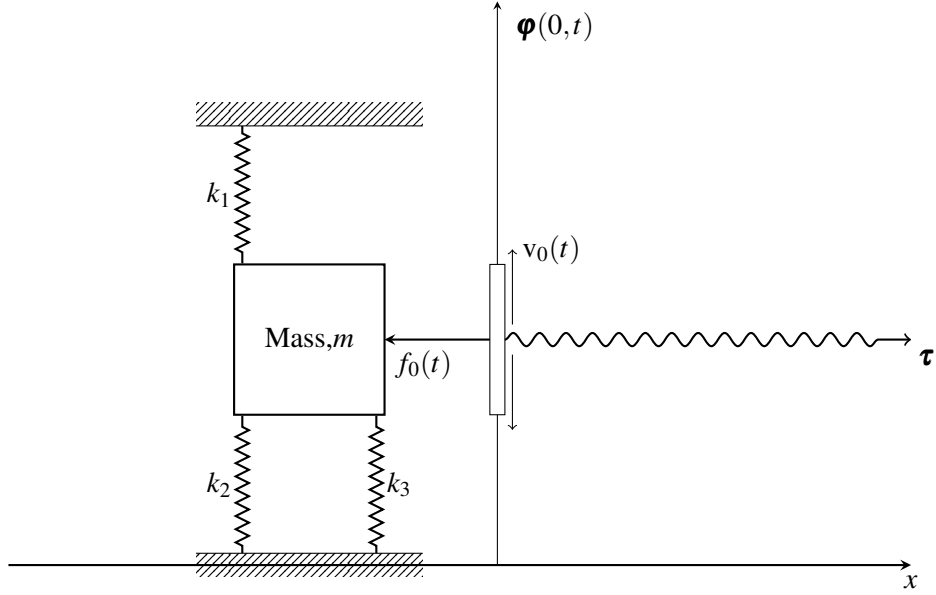
\begin{figure}[h!]
\centering
\begin{tikzpicture}[M1/.style={rectangle,draw=black,minimum size=2cm,thick}]
    \tikzstyle{spring}=[thick,decorate,decoration={zigzag,pre length=0.1cm,post length=0.1cm,segment length=5}]
    \tikzstyle{ground}=[fill,pattern=north east lines,draw=none,minimum width=0.75cm,minimum height=0.3cm]
     (-3.0, 2.5)\node  [M1] (M1) {} ;                   
    \node (wall1) [ground, minimum width=3cm,yshift=-3cm] {};
    \draw (wall1.north west) -- (wall1.north east);
    \draw [spring] (wall1.170) -- ($(M1.south east)!(wall1.170)!(M1.south west)$) node[pos=.5,left] {$k_2$};
    \draw [spring] (wall1.10) -- ($(M1.south west)!(wall1.10)!(M1.south east)$) node[pos=.5,right] {$k_3$};
    \node (wall2) [ground, minimum width=3cm,yshift=3cm] {};
    \draw (wall2.south west) -- (wall2.south east);
    \draw [spring] (wall2.190) -- ($(M1.north east)!(wall2.190)!(M1.north west)$) node[pos=.5,left] {$k_1$};
    \node  at (M1.center) {Mass,$m$};
       
\draw [-stealth, thick] (-4,-3.0)--(8,-3.0) node[anchor=north] {$x$};  
\draw [-stealth] (2.5,1)--(2.5,4.5);   
\draw (2.4,-1) -- (2.4,1) -- (2.6,1) -- (2.6,-1) --cycle;  
\draw[->] (2.7,0.2)--(2.7,1.2)node[anchor=west] {$\v_0(t)$};
\draw[->]  (2.7,-0.2)--(2.7,-1.2);

\draw (2.5,-3)--(2.5, -1);
\draw  (2.6, 4.5) node[anchor= north west, thick] {$ \varphib(0,t)$};
\draw [-stealth, thick] (2.4, 0)--(1.0,0) node[anchor=north west] {$f_0(t)$};  
\draw [->,snake=snake,line after snake=3. mm,thick]  (2.6,0) -- (8.0,0);
 \draw  (8.2, 0) node[anchor=  north] {$\taub$};
\end{tikzpicture}
 \caption{The semi-infinite string.}   
 \label{figure2}
\end{figure}

 Horace Lamb, a   professor of Fluid Mechanics  in Cambridge, in 1900, proposed a
 similar infinite-dimensional model to explain the surge of a  wave-like behavior in an infinitely long violin string.\\
 Lamb's model is  a semi--infinite string tautly stretched at tension $\tau$, connected
at the left-end point, situated at $x=0$, to a lumped {\em conservative}
mechanical load, for example composed of an arbitrary (but finite) number of
point masses and linear springs connected together as in Fig. \ref{figure2}. The first paper analyzing this system  in depth seems to be   \cite{Lewis-T-75} while in  \cite{Picci-T-92} there is  an analysis,    which will be replicated below, in the same spirit of that done before for the electrical line. 

When it becomes infinitely long, once plucked in an arbitrary point, the string does not oscillate anymore but the local deformation is instead sweeped through in both directions by deflection waves. These waves interact with the load inducing a motion which is not longer oscillatory. This motion is the object of the analysis which follows.

Let $\varphi(x,t)$ be the vertical deflection of the string at distance $x$ from
the load and let
 \begin{equation}
 \v(x,t) = \frac{\partial}{\partial t} \varphi(x,t), \quad  f(x,t) = \tau
\frac{\partial}{\partial x} \varphi(x,t)     \label{vf} 
\end{equation} 
be the  vertical components of the velocity and  of the tension of the string at
$x$.

At $x=0$ the pulling force $f_{o}(t) := f(0,t)$ acts on the mechanical load
inducing a constrained motion along the vertical axis with velocity $\v_{o}(t) :=
\v(0,t)$, related to the acting force by the mechanical impedance of the load.
Using Laplace transforms we can write  
\begin{equation}  
\hat{\v}_{o}(s) = Z_{o}(s) \hat{f_{o}}(s)     \label{IMPED}  
\end{equation}  
 where $Z_{o}(s)$ is the mechanical impedance of the load seen from the connecting
point with the string. 

The system can be described by the following state equations  
\begin{eqnarray}
\left[ \begin{array}{c}
\frac{\partial \v}{\partial t}\\ \frac{\partial
f}{\partial t}\end{array}\right] & =  & \left[ \begin{array}{cc} 0 & 1/\rho
\frac{\partial}{\partial x}\\   \tau \frac{\partial}{\partial x} & 0 \end{array}\right] 
 \left[ \begin{array}{c} \v \\f \end{array}\right]   \label{WAVE}  \\     
\dot{\xi} & = &  A\xi + b_{o}f_{o}                \label{STATE}  \\
  \v_{o} & = & c_{o}\xi, \qq \xi\in \Rbb^{2n}                    \label{VO}  
  \end{eqnarray}  
where $\rho$ is the density of the string. The last two equations (\ref{STATE})
and (\ref{VO}), can be thought as a realization of the mechanical impedance
$Z_{o}(s)$ which is   expressible as   a rational function having a minimal realization  
\begin{equation} 
 Z_{o}(s) = c_{o}(sI - A)^{-1}b_{o}         \label{ZO}  
 \end{equation}
where the matrix $A\in \Rbb^{2n\times 2n}$ has purely imaginary eigenvalues. Recall again that  since $Z_{o}(s)$ is a {\em lossless} impedance function  there
cannot be a direct feedtrough term in the equation \eqref{VO} of the realization.

 One wants to model the motion of an (observed) output variable of the system which is formed
as a linear combination of the state of the load $\xi$ with perhaps a direct
feedtrough term from the pulling force  of the string $f_{o}$, say  
\begin{equation}
 y = c\xi + df_{o} \,.     \label{OUTPUT} 
  \end{equation}
 The string equation (\ref{WAVE}) is just the wave equation written in vector
form. We assume suitable units have been chosen to insure $\tau/ \rho = 1$ so
that the speed of propagation along the string is one. The evolution of the
composite system (\ref{WAVE}),(\ref{STATE}) can  be seen as the evolution of a
conservative Hamiltonian system $\dot{z} = Fz $ with state (phase) vector
\begin{equation}
 z:= \left[\begin{array}{c} \xi \\ v \\f \end{array} \right] \label{Z} 
\end{equation} 
 taking place in the phase space ${\bf H} := {\bf R}^{2n}\oplus
L^2_{2}({\bf R_+})$. This space can be given a Hilbert space structure by
introducing the {\em energy norm} 
 \begin{equation}
 \| \left[\begin{array}{c} \xi \\ v \\f \end{array} \right] \|^2 = 1/2 \xi^{\top}\Omega
\xi + 1/2 \int^{+\infty}_{0} (\rho \v^2 + f^2) dx   \label{NORM}  
\end{equation}
where $\Omega$ is a symmetric nonnegative matrix representing the total energy
(hamiltonian) of the load, a quadratic form in the state $x$. By choosing $x$
minimally one can always guarantee  $\Omega >0$. 

The $F$ operator  in the dynamical equations
(\ref{WAVE}),(\ref{STATE}), (\ref{VO}) is skew-adjoint on its natural domain (of
smooth functions satisfying the boundary conditions (\ref{IMPED})) and generates
an  energy preserving (i.e. orthogonal) group on ${\bf H}$.

Since the string subsystem obeys the wave equation we can express the
displacement in  D'Alembert form  
\begin{equation}
 \varphi(x,t) =a(t+x) + b(t-x) \hspace{1cm}t \in \Rbb, x \geq 0 \label{AB}
\end{equation}
 where the functions $a$ and $b$ are  the {\em incoming} and {\em
outgoing} waves respectively. In the present setup it is actually only the
derivatives $a'$ and $b'$ which will enter the scattering representation of the
state vector $(\v,f)^{\top}$ as determined  by the initial data of (vertical) 
velocity and tension along the string at (say)  time zero. In fact, putting
$t=0$ we have from \eqref{vf}, (\ref{AB})  
\begin{eqnarray}
   \v_{i}(x) := \v(x,0)&=& a'(x) + b'(-x) \nonumber \\
    f_{i}(x) := f(x,0)& = & a'(x)- b'(-x) \,. \label{IC}  
    \end{eqnarray}  
This system of equations determines only
$a'(x)$ {\em for} $x \geq 0$ and $b'(x)$ {\em for} $x \leq 0$. A key point is that in this picture, the initial data $(v_{i},f_{i})^{T}$ in $L^2_{2}({\bf R}_{+})$ are {\em arbitrary} and the
restrictions $a'|_{x\geq 0}$ and $b'|_{x \leq 0}$ are therefore also arbitrary. They could be called   "free
variables" in $L^2({\bf R}_{+})$ and $L^2({\bf R}_{-})$, respectively. Later we shall model them as elementary events in some probability  space.

Now, from the same  idea of Scattering Theory as in the previous section,
we derive a mathematical descriptions  of the boundary variables
$\v_{o}(t)$ and $f_{o}(t)$ by linear{\em finite dimensional } models driven by
free $L^2$--input variables. We shall discover that these linear models  have similar structural
properties  to those of {\em stochastic} models driven by white noise processes. 

The procedure  starts with the identities  
\begin{eqnarray}  \v_{o}(t) & = & a'(t) + b'(t) \label{FSCATT}\\    
f_{o}(t) & = & a'(t) - b'(t)   \label{SCATT}  
\end{eqnarray} 
 and then uses the "steady--state" boundary condition at $x=0$, relating $\v_{o}(t)$ and
$f_{o}(t)$ specified by the load impedence (\ref{IMPED}). It is immediate to check that
(\ref{FSCATT}) and (\ref{SCATT}) are both interpretable as {\em state feedback}
laws on the mechanical load system $(A,b_{o},c_{o})$ i.e. they  can be
written   \begin{eqnarray}
 f_{o}(t) &=& -\v_{o}(t) +2a'(t) \label{FF}\\       
 f_{o}(t) &=&  \v_{o}(t) -2b'(t)         \label{BF} 
\end{eqnarray} 
 so that (\ref{FF}) and, respectively, (\ref{BF}), after
substitution of  (\ref{VO}) and combining with (\ref{OUTPUT}),  yield the
following pair of representations of the "output signal" $y$ \footnote{Of course
one can in particular obtain analogous representations also for the tension $f_{o}$ and velocity
$\v_{o}$}  
\begin{eqnarray}
 \dot{\xi}(t) & = & (A -b_{o}c_{o})\xi(t) + 2b_{o}a'(t) \label{FR}\\   
y(t) & = & [c- dc_{o}]\xi(t) + 2a'(t)   
\end{eqnarray}
 and, respectively  
\begin{eqnarray}
 \dot{\xi}(t) & = & (A + b_{o}c_{o})\xi(t) -2b_{o}b'(t)  \label{BR}\\
 y(t) & = & [c+ dc_{o}]\xi(t) - 2b'(t)  
\end{eqnarray}
Hence  any output $y$ of the system, of the type
(\ref{OUTPUT}) admits a bona fide {\em Forward--Backward pair of representations}
 in the spirit of stochastic realization of stationary
processes. In fact, the properties of the two models are exactly the same as those derived earlier  for the electrical line in Theorem \ref{MainElThm}. They can be summarized in the following statement.

\begin{thm}\label{MainThmMech}
 Assume the realization (\ref{ZO})
is minimal. Then  the representations (\ref{FR}) and (\ref{BR}), respectively,
are asymptotically stable and antistable, in fact  
\begin{equation}  
\Re\lambda(A - b_{o}c_{o}) <0,\hspace{1cm}  \sigma(A + b_{o}c_{o}) = - \sigma(A -
b_{o}c_{o})   \label{EIG} 
 \end{equation}  
 the symbol $\sigma(A)$ denoting the spectrum of the matrix $A$. 

Moreover the two representations are related by a "change of white noise input"
formula of the type $\hat{b}' = K(s)\hat{a}'$ with $K(s)$ an inner function where  $K(s)$ is precisely the scattering function associated to the boundary
condition (\ref{IMPED}), i.e.   
\begin{equation} 
 K(s) = \frac{Z_{o}(s) -1}{Z_{o}(s) + 1}\,. \label{K}  \end{equation}  
 \end{thm}

\begin{rem}In a sense, if we anticipate some     stochastic realization theory summarized later on in Chapter 3 of Part two, this representation result  can be {\em reversed} in the following sense: Start with a stationary process $y$ of assigned
rational spectral density $\Phi(s)$,  select an analytic--coanalytic pair of
spectral factors of $\Phi$, say $W(s),\, \bar W(s)$, define the relative inner scattering matrix $K(s)$ by setting
$$ K(s)= \bar W(s)^{-1}W(s)
$$
and form the lossless impedance $Z_{o}$ by solving (\ref{K}). Then the process $y$ with the given spectrum can be represented by a dynamical model which describes the connection of
the lossless load of impedance $Z_{o}$ to a lossless infinitely long string.
This stochastic finite-dimensional representation   essentially holds valid for any linear observable \eqref{OUTPUTel} of the electrical (or mechanical) load connected to an infinitely  long lossless line.
\end{rem}

Concerning the abuse of the ``white noise input" qualification which seems to have been arbitrarily attached to  the input functions  $a'$ and $b^{\prime}$, let us first remark that they are initially    determined  on positive (and respectively negative)  half lines $\{t \geq 0\}$ and $\{t \leq 0\}$ by the initial conditions
of the string,   but that they propagate in the  {\em Forward} and {\em Backward} directions of time by   the action of the forward and backward time shift operators as in \eqref{TranslICa} and \eqref{TranslICb}. In a deterministic setting they could  then be qualified as ``free input functions". We will prove in Sect. 1.5 of Part Two that their qualification as   {\em white noise processes} can  indeed be justified mathematically.

\subsection{\bf Example 3: A  particle in an infinite-dimensional heat bath}\label{sub:Ex3}

Consider Boltzmann's ideal gas model extended to an {\em infinite number of equal particles}. Assume we observe the motion of only one particular particle;  could it show an  irreversible behavior?

For simplicity we consider a one-dimensional array of equal particles of unit mass  indexed by an integer-valued label $k=0,\pm1,\pm 2, \ldots $ each performing a one dimensional motion with respect to a rest position on some one-dimensional fixed lattice, with configuration (position) variable $q_k(t)$ and momentum $p_k(t)$ at time $t$.

\begin{figure}[h!]
\centering
\begin{tikzpicture}
\begin{scope}[scale=0.9,transform shape]
\foreach \x in {-4,-3,-2,-1}  \draw(\x,0) [fill=yellow]circle(2. mm);
 \draw(0,0) [fill=red]circle(2. mm);
\foreach \x in {1, 2, 3 ,4}  \draw(\x,0) [fill=yellow]circle(2. mm);
\foreach \x in {-5, -6, -7}
\draw(\x,0) [fill=black]circle(1.pt);
\foreach \x in {5,6, 7}
\draw(\x,0) [fill=black]circle(1.pt);
\end{scope}
\end{tikzpicture}
\caption{the Brownian particle in a one-dimensional heat bath}

\end{figure}
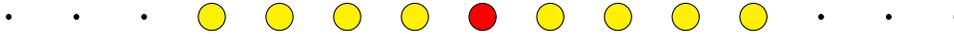

Assume the dynamics of the system is generated by a quadratic Hamiltonian
\beq \label{DiscrHamilt}
H(q,p)= \met \sum_{k}p_k^2 + \met \sum_{k,h}q_k V^2_{k,h}q_h
\eeq
where we choose  the infinite potential matrix with  a banded tridiagonal symmetric  structure
\beq \label{Toepl}
V^2:= \bmat \dots& \ldots& 0& \ldots& \ldots& \ldots\\
		a_1& a_0& a_1& 0& \ldots& \ldots\\
		\ldots& a_1& a_0& a_1&0& \ldots\\
		\ldots&	\ldots& a_1&a_0&a_1& \ldots\emat
\eeq
 with  the triplet $a_1,\, a_0,\, a_1$ proportional to  $-1,\, 2,\, -1$ so that $V^2$ is  positive definite and is bounded above by some multiple of the identity operator. The systems is conservative, assumed in steady state with the canonical variables $p(t)$ and $q(t)$ evolving in $\ell^2(\Zbb)$. We shall denote them by doubly-infinite column vectors indexed in increasing order.
 

 The gradient of the potential of the system
 $ \varphi(q) := \Frac{1}{2} q^{\top} V^2 q$  is written as the   infinite column vector $V^2 q$   whose  $k$-th component (equal to the gradient of the potential  at a distance $k$ from the origin)   is
 \beq \label{Gradk}
 \nabla\varphi(k):= [\,V^2\,]_{k} \,q= c^2\bmat -1& 2& -1\emat \bmat q_{k-1}\\q_{k}\\q_{k+1}\emat\,: =c^2 \Delta^2_k \,q, \qquad k=\pm1,\pm2,\pm3,\ldots
 \eeq
Note  that $\Delta^2_k \,q$ is  proportional to  the {\em opposite of  the discrete second spatial derivative} of the configuration variable $q$ evaluated at the point at  distance $k$ from the origin. Hence  for the force field $f(k)$ acting at distance $k$, we obtain the expression
 \beq
 f(k)= -\nabla \varphi(k):= -c^2 \Delta^2_k \,q,\qq k=1,2,3,\ldots
  \eeq
   Since $\dot {p}=f$ (or $\ddot{q}(t)= \frac{1}{m}f(t)$)  the Hamiltonian equations for the yellow particles (the ``heath bath") have the form
 \beq\label{DiscrCanEq}
  \bmat \dot{q} \\ \dot{p}\emat = \bmat 0 & 1\\ -c^2 V^2 &0\emat \bmat  q\\p\emat
  \eeq
  where   $ q$ and $p$ are infinite column vectors with components indexed by the integer $k=\pm1,\pm 2, \pm 3, \ldots$.
  
  We shall isolate the dynamics of the red particle of index $0$ from the canonical equations of the heath bath, noticing that for $k=0$    the force \eqref{Gradk}, acting on particle 0  only  depends on $q_ 0, q_1$ and  and $q_{-1}$. In other words  there is only interaction with the nearest particles of index $\pm1$, i.e.
\beq\label{BrownianP}
\begin{cases}
\dot q_0(t) & = p_0(t)\\
\dot p_0(t) &=  -2c^2\,q_{0}(t)+c^2 [q_1(t)+ q_{-1}(t)]
\end{cases}
\eeq
where we may interpret $ c^2[q_1(t)+ q_{-1}(t)]\equiv c^2 u(t)$  as a force action  from the heath bath on the $0$-th particle. 
 This system is just a one-dimensional mechanical oscillator at frequency $\omega_0=\sqrt{2c^2}$.
 Since we have normalized to unit mass   the  oscillator equation  can be written in a more familiar form as
$$
\ddot{q}_{o}(t)   =  -  \omega_0^2 q_0(t) +c^2 u(t)\q \text{or, equivalently} \qq \ddot{p}_{o}(t)   =  -   \omega_0^2 p_0(t) +c^2 \dot{u}(t)
$$
which are actually   {\em forced} oscillator equations with input force  $u(t)\equiv [q_1(t)+q_{-1}(t)]$, or $\dot{u}(t)\equiv [p_1(t)+p_{-1}(t)]$ coming from the interaction with    the closest particles of the heat bath. The system output can be chosen either as  $y=p_0$ or $y=q_0$ both leading to the transfer function
\beq \label{TrFcn}
Z_0(s)= \Frac{c^2}{s^2 +\omega_0^2}
\eeq
relating the input $u$ to $y$. The system \eqref{BrownianP} should be  considered as defining the boundary conditions at the location $k=0$ for the  equation of the heat bath. A warning  here is that the  inputs $q_{\pm1}$ are  produced by the full dynamics of the system and cannot  be understood as  ``free variables" as in the electrical line example.

The ultimate goal of our analysis should now be  to show that  in an infinite-dimensional heat bath the input variables $u(t)$ or $\dot{u}(t)$  behave (probabilistically) like a white noise input but also that the interaction of the 0-th particle with the heat bath can be described by a {\em feedback relation which changes the dynamics of the oscillator and makes it behave as a dissipative system}. This was (quite implicitly) the Ford-Kac-Mazur program \cite{Ford-K-M-65} who obviously could not try a system-theoretic justification like the one we have in mind.

 We shall  attempt an analysis similar to what we have done for the electrical line. The first step is to introduce a {\em non canonical} change of variables which leads to a symmetrical dynamics like the equation \eqref{WAVE2} for the pair $\v, i$. To this purpose, let's introduce the factorization
 \beq \label{CholFact}
 V^2= V V^{*}
 \eeq
 where we choose $V$ to be the normalized  lower triangular  (infinite)  factor matrix of $V^2$ and define
 \beq \label{ComplexBrownian}
 x(t):= V^{*} q(t),\qquad z(t):= p(t)+ j x(t)
 \eeq
 so that the Hamiltonian of the system becomes $H(z):= \met \|z\|^2$. By this transformation 
the Hamiltonian dynamics of the heat bath \eqref{DiscrCanEq} is transformed to
\beq\label{SymmetrHamilt}
\bmat \dot x(t)\\ \dot p(t)\emat = \bmat 0& V^{*} \\ -V & 0\emat  \bmat x(t)\\ p(t)\emat
\eeq
so that  (by symmetry of $V^2$)
$$
\ddot{p}(t)= -V^2 p(t),\qq \ddot{x}(t)= -V^2 x(t)
$$
which are both {\em ``semi-discrete"} wave equations for the time evolution of both infinite vectors $x$ and $p$ in $\ell^2$ of the explicit  form
  \beq \label{WaveBr}
\frac {\partial^2 x(k,t)}{\partial t^2}= - c^2 \sum_{h=k-1}^{k+1}\Delta^2_{k,h}\,x_{h}(t),\qq\frac {\partial^2 p(k,t)}{\partial t^2}= - c^2 \sum_{h=k-1}^{k+1}\Delta^2_{k,h}\,p_{h}(t)\,  \,\qq k\in \Zbb
\eeq
where $\Delta_k^2 $ is the discrete   second derivative operator  with respect to the location index  $k$, which could be associated to boundary conditions $q_0(t),p_0(t)$ at $k=0$  constrained by the dynamical equations \eqref{BrownianP}.

 This system   has evidently an  hyperbolic character and in the new (complex) phase variable the evolution must preserve the $\ell_{\Cbb}^2(\Zbb)$  norm. Hence  its solution must be  a {\em unitary} group, unitarily equivalent (see later) to  translation in time on $\ell^2$.\\
 Therefore we shall adventure to state that      both  $x$ and $p$ should  be expressible in   in terms of some {\em incoming} and {\em outgoing} waves respectively  (d'Alembert formula). Here however the discrete  structure of the space variable creates some difficulties. It is not clear how one should model a wave traveling on an infinite (discrete) lattice.
 
 { \em Formally}, we shall imagine an hypothetical fixed distance $h$ between the particles and waves traveling with some velocity $v$ so that, referring to waves observed  at  points at a distance $nh$ from the origin $x=0$ in a geometric  lattice, an incoming wave   at distance $nh$ could be denoted $\bar a(nh + vt)$. Below we shall   reserve the symbol $a(k +t)$ to denote its one-space-step increment at the location $kh$. Both signals $x(k,t)$ and $ p(k,t)$ are steady state solutions of the wave equation and therefore have a wave-like behavior of this kind.
 
 \begin{rem} Returning to the similarity with the dynamics \eqref{WAVE2} for the pair $\v, i$ in the electrical line, we may identify  the variable $x$ with voltage and $p$ with current (or conversely). This similarity and the same  structure of the Hamiltonian (up to coefficients) may lead to guess that the relation between $x(k,t)$ and $ p(k,t)$ at each finite location $k$ should be a {\em positive-real} kind of relation described by a positive-real impedance function. By a well-known transformation in network theory \cite{Anderson-V-73} this may be equivalently expressed by a {\em lossless (unitary) scattering}  transformation between two related wave functions. This is what we shall attempt to do below.
 
 Further it is worth noting that   the solution of the dynamical equations \eqref{SymmetrHamilt} is  invariant with respect to exchanging the role of $V$ and $V^{*}$ as the wave equation for both variables only depends on $V^2$.
\end{rem}

  Following a standard procedure in scattering theory, we introduce the functions $a$ and $b$
  \begin{equation} \label{AB}
a(k+t) :=\met [ x(k,t) + p(k,t)], \qq b(k-t) := \met [ x(k,t) - p(k,t)], \qq t \in {\bf R}, \;\q k \in \Zbb
\end{equation}
 which are also solutions of the semi-discrete wave equations \eqref{WaveBr} now explicitly expressed in d'Alembert form. They  evolve in time by  continuous-time translation of two $\ell^2$ sequences  $\ab:=\{a(k), k \in \Zbb\}$ and $\bb:=\{b(k), k \in \Zbb\}$ which are determined by the initial distributions of the variables $x$ and $p$ at time $t=0$, as obtained from \eqref{AB},  namely
\beq
\ab= \met (\xb + \pb),\qq \bb= \met (\xb-\pb)
\eeq
where the boldface denote sequences in $\ell^2$. The backward and forward time shift of these variables will be   called the {\em incoming} and {\em outgoing waves} respectively.\\
In order to provide some intuition for the following construction, we shall anticipate a bit  the theory exposed later in Part two, which states that,  under any absolutely continuous  invariant measure for the shift group, $\xb, \pb$ can be understood as stochastic processes depending on the discrete variable $k$. Since the Hamiltonian has the additive structure $H(\xb, \pb)= \met \|\xb\|^2 + \met\|\pb\|^2$ the two processes turn out to be   {\em independent white noises} (in this sense, stochastic {\em free variables}) with the same variance and so this will therefore  also be true for the two linear combinations $\ab$ and $\bb$ as can readily be checked. By the effect of time translation, these random sequences (i.e. discrete processes) describe a continuous time-indexed flow of elements in $\ell^2$ which, using   similar notations  to \eqref{TranslICa}, \eqref{TranslICb}   have sample values
\beq
a(k+t)= [\Sigma(t)\ab](k),\qq b(k-t)=[\Sigma^{*}(t)\bb](k),\qq t\in \Rbb\,.
\eeq
Under the invariant measure mentioned above, these wave-like signals will be  understood as  Hilbert-space valued, continuous-time   generalized stochastic processes. Their continuous-time white noise character will be unveiled later  in Sect. 1.5  of Part two.\\
 Consider in particular the incoming and outgoing waves measured at the location   $k=0$, say $a_0(t): a(0+t)$ and $b_0(-t):= b(0-t)$ which by \eqref{AB} enter in the dynamics through the relations
$$
x_0(t) + p_0(t)= 2 a_0(t),\qq x_0(t) - p_0(t)= 2 b_0(-t)
$$
which written in function of the original configuration variable $q$ take the form
\beq \label{FeedV}
V_0^{*}q(t) +p_0(t)= 2 a_0(t), \qq V_0^{*}q(t) -p_0(t) =  2 b_0(-t)\,.
\eeq
where $V_0^{*}$ denotes the row of index 0 of the factor $V^{*}$ in the Cholesky-type factorization \eqref{CholFact} of $V^2$. These will  actually turn out to be the feedback relations expressing the input $u$ as a linear combination of the state variables  plus noise. These feedback laws depend on the structure of the factor $V_0^{*}$ which we shall discuss in the following.
 \begin{thm}
The dynamics of the Brownian particle in the configuration of Fig.1 is described  by both  the  stable forward state-space model
\beq \label{stableBr}
\bmat \dot q_0(t) \\ \dot p_0(t) \emat  = \bmat 0 & 1 \\ 0&  -2c\emat \bmat   q_0(t) \\  p_0(t) \emat  +\bmat 0\\4c\emat  w(t)      
\eeq
and by the antistable backward state-space model
\beq \label{AntistableBr}
\bmat \dot q_0(t) \\ \dot p_0(t) \emat  =  \bmat 0& 1\\0 & 2c\emat \bmat   q_0(t) \\  p_0(t) \emat  +\bmat 0\\4c\emat  \bar w(t)      
\eeq
which has the opposite positive eigenvalue   $\lambda=+2c$. The input noise processes $w(t)$ and $\bar{w}(t)$ are   the incoming and out going waves $a_0(t)$ and $b_0(-t)$ at the location $k=0$.\\
The momentum $p_0(t)$ of the Brownian particle satisfies the differential equation
\beq \label{LangevinEq}
\dot{p_0}(t)= -2c p_0(t) +4c   w(t)
\eeq
 which is exactly of the same form   predicted by the classical physical analysis of Paul Langevin  \cite[eq. (3)]{Langevin-908}.
But,  $p_0$ also  satisfies the backward Langevin-type equation:
\beq \label{BackLangevinEq}
\dot{p_0}(t)= 2c p_0(t) + 4c  \bar{w}(t)
\eeq
derived from \eqref{AntistableBr}. The formal double-sided Laplace transform of the incoming and out going input signals $\hat{w}(s)$ and $\hat {\bar w}(s)$ are related by an all-pass fiter
\beq \label{InnerBrownian}
\hat{w}(s)=Q(s) \hat{\bar w}(s),\qq Q(s)=\Frac{s-2c}{s+2c}
\eeq
 where $Q(s)$ is actually inner.
 \end{thm}
\begin{proof}
We shall  present a direct  proof for the   configuration of the one-dimensional lattice  where the heat bath  particles are  on both sides of the Brownian particle. The particles are  indexed by an integer-valued  index $k$ which can take both positive and negative values on the  integer set $\Zbb$. Since the potential matrix \eqref{Toepl} (denoted by $V^2$) is   a doubly infinite symmetric Laurent matrix, its factorization can be reduced to the factorization of its {\em symbol} which we write
$$
\hat{V}^2(z)= c^2[ -z^{-1} +2 - z]= c^2 [z^{-1} -1][z-1] = \hat{V}(z)\hat{V}^{*}(z)
$$
so that $V^{*}$ is $c$ times an infinite {\em  upper triangular} matrix  whose  rows are all equal to $\bmat \ldots&0&-1& 1&0  \ldots\emat$ the main diagonal being  a sequence of  of $-1$'s and the upper diagonal a sequence of  of $+1$'s. Hence   by the ordering of the components in the vector  $q\in \ell^2(\Zbb)$ we have
\beq
[V^{*}q ]_k= c(q_{k+1}-q_k) \qq
\eeq
where  $[\cdot]_k$ stands for $k$-th component. For convenience we shall rewrite here  the state space model for the dynamics of the Brownian particle \eqref{BrownianP} 
\beq\label{BrownianP2}
\bmat \dot q_0(t) \\ \dot p_0(t)\emat= \bmat 0 &1\\ -2c^2 &0\emat \bmat q_0(t)\\ p_0(t) \emat + \bmat 0\\c^2\emat u(t)
\eeq
 where  the input force $u(t)$  at time $t$  is  identified with the  sum $q_{-1}(t)+ q_1(t)$, representing a symmetric  interaction  from both sides of the lattice.\\
After the explicit introduction of the  extended factor $V^{*}$, the relations \eqref{FeedV}  yield
$$
  c[q_1(t) -q_0(t)] +p_0(t)= 2 a_0(t), \qq  c[q_1(t)-q_0(t)]-p_0(t) =  2 b_0(-t)\,.
$$
which provide the expressions
\beq
q_1= \bmat 1 &- \frac{1}{c} \emat \bmat q_0\\p_0\emat + \frac{2}{c} a_0(t),\qq q_1= \bmat 1 &\frac{1}{c} \emat \bmat q_0\\p_0\emat +\frac{2}{c} b_0(-t)\,.
\eeq
Symmetrically, one may now repeat the derivation starting from a definition of $x$ involving the lower banded factor of $V^2$, i.e.  $x(t):= Vq(t)$ so that $x_0(t)= c(q_{-1}(t)-q_0(t))$ and one obtains
\beq
q_{-1}= \bmat 1 &-\frac{1}{c} \emat \bmat q_0\\p_0\emat + \frac{2}{c} a_0(t),\qq q_{-1}= \bmat 1 &\frac{1}{c} \emat \bmat q_0\\p_0\emat + \frac{2}{c} b_0(-t)\,.
\eeq
Combining the two relations above,  the input $u(t)= [q_1(t)+q_{-1}(t)]$ of the state space representation \eqref{BrownianP2}  has the following two {\em state feedback} expressions
\beq
u(t)= 2\bmat 1 &- \frac{1}{c} \emat \bmat q_0(t)\\p_0(t)\emat + \frac{4}{c} a_0(t),\qq u(t)= 2\bmat 1 &\frac{1}{c} \emat \bmat q_0(t)\\p_0(t)\emat  + \frac{4}{c} b_0(-t)\,
\eeq
which lead to the two forward and backward state-space models
\begin{align}
\bmat \dot q_0(t) \\\dot p_0(t)\emat  &= \bmat 0 &1 \\ 0 & -2c\emat  \bmat q_0(t)\\ p_0(t) \emat + \bmat 0\\4c\emat w(t)\\
\bmat \dot q_0(t) \\\dot p_0(t)\emat  &= \bmat 0& 1\\0 & 2c\emat \bmat q_0(t)\\ p_0(t) \emat + \bmat 0\\ 4c\emat \bar{w}(t)
\end{align}
where 
$$
w(t):=  a_0(t)\,\qq \bar{w}(t):= b_0(-t)\,.
$$
The closed loop matrices:
\beq
\Gamma:= \bmat 0 &1 \\ 0 & -2c\emat   \qq \bar {\Gamma}= \bmat 0& 1\\0 & 2c\emat\,.
\eeq 
have the nonzero eignvalues in symmetric location at  $\lambda=\mp 2c$. Eliminating $q_0$ from the forward state-space model provides the {\em forward   equation} \eqref{LangevinEq} for $p_0$
 which coincides with that predicted by the classical physical analysis by Paul Langevin \cite[eq. (3)]{Langevin-908} while  the equation for $q_0$, namely  $\dot{q}_0=p_0$ describes precisely the "wandering trajectories" of the Brownian particle described in the literature. The relation between the two forward/backward  input noises follows by combining  \eqref{LangevinEq} and \eqref{BackLangevinEq}.
 \end{proof}
 \begin{rem}\label{RemBrownian}
 The proof works nearly the same for a three-dimensional configuration and momentum variables and a three-dimensional infinite lattice.  Note that the forward state-space model \eqref{stableBr} is not asymptotically stable as the configuration variable $q_0(t)$ of the Brownian particle  is the time integral of the  "noisy" momentum $p_0(t)$. In a  stochastic description (to be addressed later)  it can be shown that its variance grows linearly in time and hence the particle evolves erratically in an unbounded manner like a Wiener process (which is a continuous time analog of a random walk).
 Dual considerations apply to the backward model.\\
Pictorially, the backward model describes the same trajectory of a smoke particle traveling back to the tip of the cigarette.
 \end{rem}

  \subsection{\bf Relations with ScatteringTheory}

 The message from the examples of this chapter is that some conservative finite-dimensional systems coupled to a conservative Hamiltonian  system with an infinite-dimensional phase space (the heat bath) can behave in a totally unexpected way, showing a stochastic-like behavior. This behavior is caused by  the action on the finite-dimensional system of a forward and a backward wave  which are the effect of  the coupling with the heat bath dynamics. These forward and backward  travelling waves are typical of the Hamiltonian hyperbolic structure  assumed for  the heat bath and in principle, besides the classical wave equation,   may   be describing  other linear infinite-Hamiltonian systems like that introduced for the Brownian particle in Example 3.
 
Geometrically, one can see  the overall Hilbert space of the composite system as the vector sum at time $t=0$
$$
 \Hb= \Sb \vee \bar {\Sb}, \qq \Xb = \Sb \cap \bar {\Sb}
$$
where $\Sb$ and $\bar{\Sb}$ are incoming and outgoing subspaces, respectively left- and right-invariant for the Hamiltonian evolution operator $\Phi(t)$, which are generated by the past and future  (at time zero) of the "noise signals" $w(t)$ and $\bar{w}(t)$. In standard Hilbert space notation $\Sb= \Hb^-(w)$ (the past space of $w$ at time 0) and $\bar {\Sb} =\Hb^+(\bar w)$ (the future space of $\bar w$ at time 0). The subspace $\Xb$ is the finite dimensional subspace of $\Hb$ linearly generated by the components of the state vector, say at time zero, $\xi(0)$, which in the specific examples can be easily  seen as being simultaneously linear functions of either  the infinite past of  $w$ or of the infinite future  of  $\bar w$. See the equations \eqref{Langevin} and \eqref {LangevinBar}. \\
There is a connection of this setting with the formalism of Scattering Theory  \cite{Lax-P-67} and one may interpret the finite dimensional system as a ``scattering object" hit by the two waves. The novelty  here is that the result of the  interaction of a lossless finite-dimensional  system with an infinite dimensional (necessarily lossless) surrounding  environment produces   dissipation and stochastic-like behavior. The essential finding in the discussion of the examples of this chapter is that any linear observable {\em output} signal admits such a representation. \\
An appropriate notion in this respect  is the following.
\begin{defn}
An output signal $\yb$ of a linear Hamiltonian system   with a possibly  infinite-dimensional  Hilbert  phase space $\Hb$,  is just a linear,  possibly vector-valued,  linear functional $\cb$ on $\Hb$, whose  image $ y(t)= \cb(\pb(t),\qb(t)), t\in \Rbb$ is a function of time called the {\em observable  output} of the system. 
\end{defn}
In most cases of physical interest, like e.g. the infinite electrical line or the infinite particle system,  it   only makes sense to consider a {\em finite number} $m$ of scalar observables which, in the present setting, we shall  write in  vector notation as:   
$$
z\,\to\, \cb(z):= \bmat \langle \cb_1,z\rangle & \ldots&\langle \cb_m,z\rangle \emat^{\top}, \quad z\in \Hb.
$$ 
For example, for the  Brownian  particle system using the complex notation of \eqref{ComplexBrownian}, $z= p+jx$, the representative of the output functional can be written as $\cb= \bmat \ldots&0&1&0& \ldots \emat \in \ell_{\Cbb}^2$ with the   $1$ in the zeroth position.

\nd{\em Non observable states} belong to  the nullspace of the observation functional  for all times $t$. They  are irrelevant to the observer so that we can and shall only consider   the dynamics restricted to the closed  invariant subspace 
\beq \label{FinMult}
\Hb_0:=\overline{\Span}_{\{t\in \Rbb\} } \{ \Phi(t) \cb_k\,; k=1,\ldots,m\}.
\eeq
and assume the equality $\Hb\equiv \Hb_0$ (a sort of infinite-time observability) which implies that  $\Phi(t)$ and its generator  have {\bf finite multiplicity} $\leq m$. This is a key condition which leads to a spectral representation of  these operators of  finite dimension.

 The concept of multiplicity is discussed in \cite[vol II, p. 913-16 ]{Dunford-Schwartz-64} or in Halmos Hilbert space book \cite{Halmos-57}. It is actually a module-theoretic concept, see   \cite[p. 91-92]{LPBook}.

\subsection{Generalization}

We are now led to address the following question. When can a linear observable of a linear infinite-dimensional Hamiltonian  system defined on a Hilbert phase space  $\Hb$, satisfying the  above observability condition, admit a finite-dimensional, dissipative  and {\em stochastic-like}\footnote{That is of the generalized Langevin structure defined above.}  dynamical representations?
 
 In the  examples of  the previous subsections  the overall  linear Hamiltonian system  was already defined as a direct sum of an infinite-dimensional  heat bath plus a finite-dimensional lossless part.  This question is obviously a bit more general as the finite dimensional lossless subsystem has now to found. It has been also called  an {\em aggregation problem}.
The problem has been answered in \cite{Picci-89} in a completely deterministic context. This paper also contains    a formal definition of forward and backward representations in a deterministic context. However since  in the stochastic context, the problem turns   actually out to be just isomorphic to  the fundamental problem of {\em stochastic realization theory} we  shall discuss its solution and the various   connections with irreverisbility in detail in Chapter 3 of Part 2.


\begin{thebibliography}{99}

\bibitem{Anderson-V-73} B. D.O. Anderson and S. Vongpanitlerd, {\em Network Analysis and Synthesis : a modern Systems Theory approach}, Prentice Hall (1973) reprinted by Dover in 2006.

\bibitem{Bernstein-011} J. Bernstein, An entropic story, in 
 	arXiv:1109.3448 [physics.hist-ph], (2011)

\bibitem{Boltzmann-895} L. Boltzmann, On Certain Questions of the Theory of Gases, {\em Nature}, vol 51, 413 (1895)

\bibitem{Boltzmann-897}  L. Boltzmann, On Zermelo's Paper "On the Mechanical Explanation of Irreversible Processes"
[Originally published under the title: "Zu Hrn. Zermelo's Abhandlung Über die mechanische Erklärung irreversibler Vorgange", {\em Annalen der Physik}, vol.  60, pp. 392–8 (1897)].


\bibitem{Clausius-865} Clausius, R.  Ueber verscheidene f\"ur die anwendug bequeme formen der hauptgleichungen der mechanischen w\"armtheorie. {\em Annalen der Physik und 
Cehemie}, 125(7): 353–400, 1865.

\bibitem{Devaney-89} R.L. Devaney: {\em An introduction to chaotic dynamical
systems}, second edition, Addison Wesley, New York, 1989.

\bibitem{Davey-11} 
Davey K. Thermodynamic Entropy and Its Relation to Probability in Classical Mechanics. {\em Philosophy of Science};78(5):955-975, 2011. doi:10.1086/662559


\bibitem{Dunford-Schwartz-64} N. Dunford and J. T. Schwartz {\em Linear Operators Part II}, Wiley Interscience 1963.

\bibitem{Fuhrmann-81} P.A. Fuhrmann {\em Linear Systems and Operators in Hilbert Spaces},
Mc Graw-Hill, 1980.

 

\bibitem{Ford-K-M-65} G.W.Ford, M.Kac, P.Mazur: Statistical mechanics of assemblies of
coupled oscillators, {\em Journal Math. Phys.} {\bf 6}, pp 504-515, 1965.

\bibitem{Foster-24} Foster, R M, "A reactance theorem",{\em Bell System Technical Journal}, Vol. 3, pp. 259–267, 1924.


\bibitem{Gelfand-V-64} I.M Gel'fand and N. Ya Vilenkin, {\em Generalized Functions, Vol 4 : Applications of Harmonic Analysis},
Academic Press, 1964.

\bibitem{Godefroy-S-91} G. Godefroy and J.H. Shapiro: Operators with dense, invariant, 
cyclic vector manifolds, {\em J. Funct. Analysis}, 1991.

\bibitem{Goldstein-80} H. Goldstein {\em Classical Mechanics, 2nd edition}, Addison Wesley, 1980.

\bibitem{Grosse-M-011}  Karl-G. Grosse-Erdmann , Alfred Peris Manguillot, {\em Linear Chaos}, Springer 2011.

\bibitem{Halmos-57} P. Halmos {\em Introduction to Hilbert Space and the Theory of Spectral Multiplicity, second ed.} Chelsea, N.York 1957.

\bibitem{Hida-80} T. Hida: {\em Brownian Motion} Springer, New York,~1980.

\bibitem{Herrero-92} D. Herrero: Hypercyclic operators and chaos, {\em Journal of Operator Theory}, Vol. 28, No. 1 (Summer 1992), pp. 93-103.
\bibitem{Jaynes-65} Jaynes, E. T.  “Gibbs vs Boltzmann Entropies.” {\em American Journal of Physics}  33:391–98, 1965..

\bibitem{Kuo-75} H.H. Kuo: {\em Gaussian measures in Banach spaces}, Springer
Verlag, New York, 1975.

\bibitem{Langevin-908} P.  Langevin, Sur la th\`eorie du mouvement Brownien. {\em Comptes rendus de l'acad\`emie des
sciences de Paris}, pp. 530-533  (1908).

\bibitem{Lax-P-67} P.D. Lax and R.S. Phillips: {\em Scattering Theory}, Ac. Press, New
York, 1967

\bibitem{Lebowitz-93-99}
J.L. Lebowitz, Boltzmann's Entropy and Time's Arrow, {\em Physics Today}, 46, 32–38(1993); see also letters to the editor and response in {\em Physics Today}, 47, 113-116 (1994);{c)} Microscopic Origins of Irreversible Macroscopic Behavior, {\em Physica A}, 263, 516–527, (1999); 

\bibitem{Lebowitz-99b}	J.L. Lebowitz,  A Century of Statistical Mechanics: A Selective Review of Two Central Issues, {\em Reviews of Modern Physics}, 71, 346–357, 1999; {e)} From Time-symmetric Microscopic Dynamics to Time-asymmetric Macroscopic Behavior: An Overview,  in {\em European Mathematical Publishing House}, ESI Lecture Notes in Mathematics and Physics. 

\bibitem{Lebowitz-08} Joel L. Lebowitz , Time's Arrow and Boltzmann's Entropy, {\em Scholarpedia}, 3(4):3448.	doi:10.4249/scholarpedia (2008). 



\bibitem{Lewis-T-75} J.T.Lewis and L.C.Thomas: How to make a heat bath, in {\em
Functional Integration and its Applications}, A.M.Arthurs ed, p. 97-123,  Oxford, Clarendon,
1975.

\bibitem{Lewis-M-84} J.T.Lewis and H. Maassen: Hamiltonian models of classical and
quantum stochastic processes, in {\em Spriger Lect. Notes in Math.} n.{\bf
1055}, pp 245-276, 1984.


\bibitem{Lindquist-P-79} A Lindquist and G. Picci, On the Stochastic Realization Problem {\em SIAM Journal on Control and Optimiz.}, {\bf
17}, pp. 365-389, (1979). doi.org/10.1137/03170.


\bibitem{Lindquist-P-85} A Lindquist and G. Picci, Realization Theory  for multivariate
stationary Gaussian processes {\em SIAM Journal on Control and Optimiz.}, {\bf
23}, pp1-50, 1985.

\bibitem{Lindquist-P-85b}  A Lindquist and G. Picci, Forward and Backward semimartingale representations of stationary increments processes, {\em Stochastics}, Vol. {\bf 15}, No. 5 ,pp. 1-50, 1985. 

\bibitem{Lindquist-P-91} A Lindquist and G. Picci, A geometric approach to modeling and
estimation of linear stochastic systems {\em Journal of Math. Syst. Estimation
and Control}, {\bf 1},pp. 241-333,1991.

\bibitem{LPBook} A Lindquist and G. Picci,{\em  Linear Stochastic Systems: A geometric approach to modeling 
estimation and identificationof linear stochastic system}, Springer Verlag 2015.

\bibitem{Massen-90} H. Maassen : Hamiltonian models of classical and quantum stochastic
processes  in {\em Realization and modeling in System Theory}, pp 505-511,
Birkhauser Boston, 1990.

\bibitem{MAA-thesis} H. Maassen, On a class of quantum Langevin equations and
the question of approach to equilibrium. {\em Thesis}, University of Groningen,
1982.


\bibitem{Nelson-58} E. Nelson, The adjoint Markoff process, {\em  Duke Mathematical Journal}, vol. 25 (1958), 671–690.

\bibitem{Nelson-67}  E.  Nelson, {\em  Dynamical Theories of Brownian Motion}. Princeton University Press (1967).\\
url:https://web.math.princeton.edu/~nelson/books/bmotion.pdf 

\bibitem{Picci-86} G. Picci: Application of stochastic realization theory to a
fundamental problem of statistical physics, in {\em Modeling, Identification and
Robust Control}, C.I.Byrnes and A.Lindquist eds., North Holland,1986.

\bibitem{Picci-88} 	G. Picci: Hamiltonian representation of stationary processes
	in {\em Contributions to Operator Theory and its applications}   I. Gohberg,
J.W. Helton, L. Rodman eds.,pp.193-215,  Birkhauser, Boston, 1988.


\bibitem{Picci-88bis} 	G. Picci:
	Stochastic aggregation, in {\em Linear Circuits, Systems and Signal Processing}  C.I.  Byrnes, C.
Martin, R. Saeks eds.  (Proceedings of  the 8th International Symposium on the
Mathematical Theory  of Networks and Systems, Phoenix, Arizona), North Holland, pp.493-501, 1988.

\bibitem{Picci-89} G. Picci, Aggregation of linear systems in a completely
deterministic framework, in {\em Three decades of Mathematical Systems Theory: A
Collection of Surveys at the occasion of the Fiftieth Birthday of Jan C.
Willems}, H. Neijmeijer, J.M. Schumacher eds. {\em Springer Lect. Notes in
Control and Information Sciences}, {\bf 135}, pp 358-381, 1989.


\bibitem{Picci-T-90} G. Picci and T.J. Taylor, Stochastic aggregation of linear
Hamiltonian systems with microcanonical distribution, in {\em Realization and
modeling in System Theory}, pp 513-520, Birkhauser Boston, 1990.


\bibitem{Picci-91}  G. Picci: Markovian representation of linear Hamiltonian systems, in
{\em  Probabilistic Methods in Mathematical Physics},
Western Periodicals, Singapore, 1991.


\bibitem{Picci-92}  G.Picci: Stochastic model reduction by aggregation in {\em
Systems Models and Feedback: Theory and Applicatons}, pp. 169-177, Birkhauser
Boston, 1992.


\bibitem{Picci-T-92} G. Picci and T. J. Taylor: Generation of Gaussian Processes and Linear Chaos, in
{\em Proceedings of the 31st Conf. on Decision and Control} IEEE Press, Tucson, AZ, 1992, pp 2125-2131. 




\bibitem{Prigogine-80} I. Prigogine, {\em From Being to Becoming}, W. H. Freeman, San Francisco, 1980.

\bibitem{Protter-87}  P.   Protter: Reversing Gaussian semimartingales without Gauss. {\em Stochastics} 20; 39-49, (1987).

\bibitem{Rozanov-67} Y.A. Rozanov, {\em Stationary Random Processes}, Holden Days, S.
Francisco, 1967.

\bibitem{Schrodinger-51} E. Schr\"odinger, Irreversibility,
{\em Proceedings of the Royal Irish Academy} vol.53, 1951.

\bibitem{Steckline-83} V.S. Steckline, Zermelo, Boltzmann and the recurrence paradox. {\em Am. J. Physics}, vol. 51 (1983) 894-897.

\bibitem{Willems-72} Willems, J.C. Dissipative dynamical systems part I: General theory. {\em Arch. Rational Mech. Anal.} 45, 321–351 (1972). https://doi.org/10.1007/BF00276493

\bibitem{Willems-03} Jan C. Willems:  Reflections on Fourteen Cryptic Issues Concerning the Nature of Statistical Inference, {\em International Statistical Review}, Vol. 71, Issue 2, (Aug 2003) , pp. 277-318.

\bibitem{Willems-86}  Jan C. Willems: From time series to linear system—Part I. Finite dimensional linear time invariant systems, {\em Automatica} 1986, p. 567.

\bibitem{Zermelo-896} E. Zermelo, Ann. Physik 57, 485 (1896). 

\bibitem{Zermelo-896b} E. Zermelo, Ann. Physik 59, 793 (1896).
 


\end{thebibliography}
\end{document}